\documentclass[11pt, a4paper]{amsart}
\usepackage{mathbbol}
\usepackage{amsmath, amsthm, amssymb}
\usepackage{txfonts, mathrsfs}
\usepackage{psfrag}
\usepackage{graphicx}

\usepackage{color}

\newtheorem*{thm*}{Theorem}

\newtheorem{thm}{Theorem}[section]
\newcommand{\bt}{\begin{thm}}
\newcommand{\et}{\end{thm}}

\newtheorem{cor}[thm]{Corollary}
\newcommand{\bc}{\begin{cor}}
\newcommand{\ec}{\end{cor}}

\newtheorem{lem}[thm]{Lemma}
\newcommand{\bl}{\begin{lem}}
\newcommand{\el}{\end{lem}}

\newtheorem{prop}[thm]{Proposition}
\newcommand{\bp}{\begin{prop}}
\newcommand{\ep}{\end{prop}}

\newtheorem{defn}[thm]{Definition}
\newcommand{\bd}{\begin{defn}}      
\newcommand{\ed}{\end{defn}}

\newtheorem{rmrk}[thm]{Remark}
\newcommand{\br}{\begin{rmrk}}
\newcommand{\er}{\end{rmrk}}

\newtheorem{quest}[thm]{Question}
\newcommand{\bq}{\begin{quest}}
\newcommand{\eq}{\end{quest}}

\newcommand{\N}{\mathbb{N}}

\newcommand{\R}{\mathbb{R}}

\newdimen\vintkern\vintkern12pt
\def\vint{-\kern-\vintkern\int}

\newcommand{\hm}{{\mathcal H}}

\newcommand{\dist}{\operatorname{dist}}
\newcommand{\diam}{\operatorname{diam}}
\newcommand{\trace}{\operatorname{tr}}
\newcommand{\length}{\ell}
\newcommand{\Area}{\operatorname{Area}}

\newcommand{\md}{\operatorname{md}}

\newcommand{\deltalip}{\delta^{\operatorname{Lip}}}
\newcommand{\fillarealip}{{\operatorname{Fill\,Area}^{\operatorname{Lip}}}}
\newcommand{\Fillarea}{\operatorname{FillArea}}
\newcommand{\Fillarealip}{\operatorname{FillArea^{\operatorname{Lip}}}}

\newcommand{\jac}{{\mathbf J}}
\newcommand{\ap}{\operatorname{ap}}
\newcommand{\apmd}{\ap\md}

\def\XXint#1#2#3{{\setbox0=\hbox{$#1{#2#3}{\int}$}
		\vcenter{\hbox{$#2#3$}}\kern-.5\wd0}}

\begin{document}

\title[]{Spaces with almost Euclidean Dehn function}

\author[S. Wenger]{Stefan Wenger}

\address
  {Department of Mathematics\\ University of Fribourg\\ Chemin du Mus\'ee 23\\ 1700 Fribourg, Switzerland}
\email{stefan.wenger@unifr.ch}

\keywords{Dehn function, isoperimetric inequality, non-positive curvature, Gromov hyperbolicity, ultralimits, asymptotic cones, Plateau problem, Sobolev maps}
\subjclass{53C23, 20F65, 49Q05}

\date{\today}

\thanks{Research supported by Swiss National Science Foundation Grants 153599 and 165848.}

\begin{abstract}
We prove that any proper, geodesic metric space whose Dehn function grows asymptotically like the Euclidean one has asymptotic cones which are non-positively curved in the sense of Alexandrov, thus are ${\rm CAT}(0)$. This is new already in the setting of Riemannian manifolds and establishes in particular the borderline case of a result about the sharp isoperimetric constant which implies Gromov hyperbolicity. Our result moreover provides a large scale analog of a recent result of Lytchak and the author which characterizes proper ${\rm CAT}(0)$  in terms of the growth of the Dehn function at all scales. We finally obtain a generalization of this result of Lytchak and the author. Namely, we show that if the Dehn function of a proper, geodesic metric space is sufficiently close to the Euclidean Dehn function up to some scale then the space is not far (in a suitable sense) from being ${\rm CAT}(0)$ up to that scale.
\end{abstract}

\maketitle

\section{Introduction and statement of main results}

The Dehn function, also known as the filling area or isoperimetric function, measures how much area is needed to fill closed curves of a given length in a space or a group by disc-type surfaces. It is a basic invariant in analysis and geometry and plays an important role particularly in large scale geometry and geometric group theory. It is a quasi-isometry invariant of a space and is connected to the complexity of the word problem in a group.

The aim of the present article is to study the geometry of spaces whose Dehn functions are sufficiently close to the Euclidean one in a suitable sense. Our study is partly motivated by the fact that the (large scale) geometry of spaces with quadratic Dehn function is not yet well understood. In general, only few properties of such spaces are known. In contrast, spaces with linear Dehn function at large scales are well understood. They are exactly the Gromov hyperbolic spaces by an important theorem of Gromov \cite{Gromov-hypgroups}. By the same theorem of Gromov, which has inspired alternative proofs in \cite{Ols91}, \cite{Bow95}, \cite{Pap95}, there are no spaces with Dehn function of super-linear sub-quadratic growth. 

In order to state our results, let $(X,d)$ be a complete metric space. Define the Lipschitz filling area of a Lipschitz curve $c\colon S^1\to X$ by
\begin{equation*}
 \fillarealip(c):= \inf\left\{\Area(v) : \text{$v\colon \overline{D}\to X$ is Lipschitz, $v|_{S^1}=c$}\right\},
\end{equation*}
where $\overline{D}$ denotes the closed unit disc in $\R^2$. See Section~\ref{sec:Sobolev-maps} for the definition of the parametrized Hausdorff area $\Area(v)$. Here, we only mention that if $v$ is injective then $\Area(v)$ equals the Hausdorff $2$-measure of the image of $v$; moreover, if $X$ is a Riemannian manifold then $\Area(v)$ coincides with the parametrized area obtained by integrating the Jacobian of the derivative of $v$.
The Lipschitz Dehn function of $X$ is the function
$$\deltalip_X(r)=\sup\left\{\fillarealip(c) : \text{$c\colon S^1\to X$ is Lipschitz, $\length(c)\leq r$}\right\}$$ for every $r\geq 0$, where $\length(c)$ denotes the length of $c$.

The following result, proved by the author in \cite{Wen08-sharp}, generalizes and strengthens Gromov's result \cite{Gromov-hypgroups} mentioned above: if a complete, geodesic metric space $X$ satisfies
 \begin{equation*}
  \limsup_{r\to\infty} \frac{\deltalip_X(r)}{r^2} <\frac{1}{4\pi}
 \end{equation*}
  then $X$ is Gromov hyperbolic. The constant $\frac{1}{4\pi}$ is optimal in view of the Euclidean plane. In the present paper we establish the borderline case of this result in the setting of proper metric spaces by proving:

\bt\label{thm:asymp-CAT-intro}
If a proper, geodesic metric space $X$ satisfies
\begin{equation}\label{eq:asymp-Eucl-isop}
\limsup_{r\to\infty} \frac{\deltalip_X(r)}{r^2} \leq \frac{1}{4\pi}
\end{equation}
then every asymptotic cone of $X$ is a ${\rm CAT}(0)$-space.
\et

Recall that a geodesic metric space $X$ is ${\rm CAT}(0)$ if every geodesic triangle in $X$ is at least as thin as a comparison triangle in the Euclidean plane, see e.g.~\cite{BrH99} for the theory of ${\rm CAT}(0)$-spaces. A metric space is proper if all of its closed bounded subsets are compact. Having only ${\rm CAT}(0)$ asymptotic cones is equivalent to geodesic triangles in $X$ satisfying the ${\rm CAT}(0)$ thinness condition up to an additive error which is sublinear in the diameter of the triangle, see \cite{Kar-Aditi-2011}. 

Theorem~\ref{thm:asymp-CAT-intro} is new even when $X$ is a Riemannian manifold and the constant $\frac{1}{4\pi}$ is optimal, see Section~\ref{sec:asymp-CAT}. The converse to the theorem does not hold, not even in the class of geodesic metric spaces biLipschitz homeomorphic to $\R^2$, as shows our next result.

\bt\label{thm:no-converse}
 There exists a geodesic metric space $X$ biLipschitz homeomorphic to $\R^2$ which satisfies
 \begin{equation}\label{eq:no-converse-ineq}
  \liminf_{r\to\infty}\frac{\deltalip_X(r)}{r^2}>\frac{1}{4\pi}
 \end{equation}
 and whose unique asymptotic cone is (isometric to) the Euclidean plane.
\et

Theorem~\ref{thm:asymp-CAT-intro} can also be viewed as a large scale analog of (one direction in) the main result of Lytchak and the author in \cite{LW-isoperimetric} which characterizes proper ${\rm CAT}(0)$-spaces in terms of the growth of the Dehn function at all scales. The proof of Theorem~\ref{thm:asymp-CAT-intro} relies on the results and proofs in \cite{LW-isoperimetric}. However, it is conceivable that arguments similar to the ones developed in the present article can be used to prove the above mentioned characterization of proper ${\rm CAT}(0)$-spaces in \cite{LW-isoperimetric} even without the condition on properness.

Theorem~\ref{thm:asymp-CAT-intro} does not only hold for asymptotic cones but also has an analog for ultralimits of sequences of proper, geodesic metric spaces $X_n$ whose Lipschitz Dehn functions satisfy
\begin{equation}\label{eq:intro-general-growth-Xn}
\deltalip_{X_n}(r) \leq \frac{1+\varepsilon_n}{4\pi}\cdot r^2 + \varepsilon_n
\end{equation}
 for all $r\in (0,r_0)$ and $n\in\N$, where $(\varepsilon_n)$ is a sequence of positive real numbers tending to zero and $r_0$ is positive and possibly infinite. We refer to Theorem~\ref{thm:main-CAT} for this analog and state here the following consequence which generalizes the result in \cite{LW-isoperimetric} mentioned above. More consequences of Theorem~\ref{thm:main-CAT} will be discussed in Section~\ref{sec:asymp-CAT}.

\bt\label{thm:intro-bigger-const-CAT}
 For every $\nu\in(0,1)$ there exists $C>\frac{1}{4\pi}$ with the following property. Let $r_0>0$ and let $X$ be a proper, geodesic metric space satisfying $$\deltalip_X(r) \leq C\cdot r^2$$ for all $r\in (0,r_0)$. Then every geodesic triangle in $X$ of perimeter $s<(1-\nu)r_0$ is, up to an additive error of at most $\nu s$, no thicker than its comparison triangle in $\R^2$. 
\et

The perimeter of a geodesic triangle is the sum of the lengths of its sides. Geodesic triangles which are, up to an additive error of at most $\nu'$, no thicker than their comparison triangles in $\R^2$ were termed ${\rm CAT}(0, \nu')$ in \cite{DeG08}. See also Section~\ref{sec:coarse-CAT} below. Aspects of this notion of coarse non-positive curvature and variants thereof were studied for example in \cite{Gro01}, \cite{DeG08}, \cite{BuF12}, \cite{Kar-Aditi-2011}, \cite{BH13}. A closely related notion is that of bolic spaces introduced in \cite{KS94}.
Proper, geodesic metric spaces satisfying $\deltalip_X(r) \leq C\cdot r^2$ for some (arbitrarily large) $C$ and all $r$ are known to be $\alpha$-H\"older $1$-connected for every $\alpha<1$ by \cite{LWY16}. It is not known whether they are actually Lipschitz $1$-connected. In view of the theorem above it would thus be interesting to study geometric properties of spaces satisfying the conclusion of Theorem~\ref{thm:intro-bigger-const-CAT} for sufficiently small $\nu$.

As already mentioned, the proof of Theorem~\ref{thm:asymp-CAT-intro} (and also of its generalization to ultralimits mentioned above) relies on the proof of the characterization in \cite{LW-isoperimetric} mentioned above of proper ${\rm CAT}(0)$-spaces in terms of the growth of their Dehn function. One of the main new ingredients established in the present paper is a solution of Plateau's problem in ultralimits of sequences of proper, geodesic metric spaces $X_n$ satisfying \eqref{eq:intro-general-growth-Xn}. We refer to Theorem~\ref{thm:energy-min-ultralimits} for the precise result, which together with the techniques and proofs from \cite{LW-isoperimetric} will yield our main result. In this introduction we only state a version of this theorem for asymptotic cones. Recall that in the context of metric spaces, the classical problem of Plateau of finding area minimizing discs with prescribed boundary was solved for proper metric spaces in \cite{LW15-Plateau} and for a certain class of locally non-compact metric spaces in \cite{GW17-Plateau}. Ultralimits and asymptotic cones of the spaces we are interested in typically fail to be proper and also do not fall into the class of spaces considered in \cite{GW17-Plateau}.

\bt\label{thm:Plateau-asymp-intro}
Let $X$ be a proper, geodesic metric space satisfying $$\limsup_{r\to\infty}\frac{\deltalip_X(r)}{r^2}\leq \frac{1}{4\pi},$$ and let $X_\omega$ be an asymptotic cone of $X$. Then every rectifiable Jordan curve in $X_\omega$ bounds an area minimizing disc which is moreover parametrized conformally.
\et

We refer to Section~\ref{sec:energy-mins} for the definitions relevant for the theorem. We do not know whether an analog of the theorem remains true when the constant $\frac{1}{4\pi}$ is replaced by a strictly bigger constant.

We end this introduction with a brief outline of the proof of Theorem~\ref{thm:Plateau-asymp-intro}. The proof of its analog for ultralimits is almost the same. In a first step we construct a candidate for an area minimizer with prescribed rectifiable Jordan boundary $\Gamma$ in $X_\omega$. Roughly speaking, this candidate comes as the ultralimit of a sequence of area minimizers in $X$ whose boundaries form a sequence of Lipschitz curves in $X$ approximating $\Gamma$ in $X_\omega$. This step relies on the results in \cite{LW15-Plateau} and \cite{LW-intrinsic} on the existence, regularity and equi-compactness of area minimizers in proper metric spaces. In a second step we show that the so found candidate is indeed an area minimizer. This is the more difficult part. The main problem is that, in general, it is not clear whether a given filling of $\Gamma$ in $X_\omega$ comes from a sequence of fillings (of suitable curves) in $X$ with almost the same area. We solve this problem by suitably discretizing a given filling of $\Gamma$ in $X_\omega$. More precisely, we show in Theorem~\ref{thm:discretization} that for every $\varepsilon>0$ and every sufficiently large $\lambda>0$ there exists a finite family $\Gamma_1,\dots, \Gamma_k$ of pairwise separated smooth Jordan curves in the open disc $D$ and a $\lambda$-Lipschitz map $\varphi\colon S^1\cup\Gamma_1\cup\dots\cup \Gamma_k\to X_\omega$ with the following properties: the restriction of $\varphi$ to $S^1$ parametrizes $\Gamma$ and $$\sum_{i=1}^k \frac{1}{4\pi}\cdot \length(\varphi|_{\Gamma_i})^2\leq \Fillarea(\Gamma) + \varepsilon,$$ where $\Fillarea(\Gamma)$ is the Sobolev filling area of $\Gamma$. Moreover, the pairwise disjoint Jordan domains enclosed by $\Gamma_i$ cover all of $D$ except a set of measure at most $\frac{\varepsilon}{\lambda^2}$. Since the Dehn function of $X$ is nearly the optimal Euclidean one on large scales and since Lipschitz curves in $X_\omega$ can be well approximated by Lipschitz curves in $X$, this discretization procedure together with Lipschitz extension arguments yield competitors in $X$ with area arbitrarily close to $\Fillarea(\Gamma)$. The existence of such competitors will then imply that our candidate minimizes area among all Sobolev discs with boundary $\Gamma$.

\medskip

The structure of the paper is as follows. In Section~\ref{sec:prelims} we collect basic definitions regarding ultralimits, asymptotic cones, the coarse ${\rm CAT}(0)$-conditon, and Sobolev mappings with values in a metric space. We furthermore recall the definition of the Sobolev Dehn function from \cite{LWY16} and a result about the existence of suitable thickenings of a metric space. In Section~\ref{sec:ET} we show that a complete, geodesic metric space whose Sobolev Dehn function is bounded by the Euclidean one must have the so-called property (ET) of Euclidean tangents. This extends a corresponding result in \cite{LW-isoperimetric} from the setting of proper metric spaces to that of general ones and is needed in the sequel. The purpose of Section~\ref{sec:discretization} is to prove that fillings in a geodesic metric space with property (ET) can be discretized in the way described above. This is used in Section 5 to construct competitors in $X_n$ starting from a  filling of a curve in the ultralimit $X_\omega$ of a sequence $(X_n)$ of spaces $X_n$ with almost Euclidean Dehn function. In Section~\ref{sec:energy-mins} we use the results from the previous sections to prove that every rectifiable Jordan curve in $X_\omega$ bounds a Sobolev disc of minimal energy and that every energy minimizer is an area minimizer. In Section~\ref{sec:asymp-CAT} we state and prove our main result, Theorem~\ref{thm:main-CAT}, which in particular implies Theorems~\ref{thm:asymp-CAT-intro} and \ref{thm:intro-bigger-const-CAT}. We furthermore prove the remaining results and discuss further consequences of Theorem~\ref{thm:main-CAT}.

\medskip 

{\bf Acknowledgements:} I wish to thank the anonymous referee for very useful comments which have led to several improvements.

\section{Preliminaries}\label{sec:prelims}

\subsection{Basic notation and definitions}

The Euclidean norm of a vector $v\in\R^n$ will be denoted by $|v|$. The open unit disc in $\R^2$ is denoted by $$D:=\{z\in\R^2 : |z|<1\},$$ its closure by $\overline{D}$ and its boundary by $S^1$.

Let $(X,d)$ be a metric space. The open ball in $X$ of radius $r>0$ and center $x\in X$ is denoted by $B(x,r):= \{x'\in X: d(x,x')<r\}$, the closed ball by $\bar{B}(x,r):= \{x'\in X: d(x,x')\leq r\}$. The space $X$ is proper if every closed ball of finite radius in $X$ is compact. Given subsets $A\subset B\subset X$ and $\nu>0$ we say that $A$ is $\nu$-dense in $B$ if for every $y\in B$ there exists $x\in A$ with $d(x,y)\leq \nu$.

 A curve in $X$ is a continuous map $c\colon I\to X$, where $I$ is an interval or $S^1$. If $I$ is an interval then the length of $c$ is defined by $$\length(c):= \sup\left\{\sum_{i=0}^{k-1}d(c(t_i), c(t_{i+1})) : \text{ $t_i\in I$ and $t_0<t_1<\dots < t_k$}\right\}$$ and an analogous definition applies in the case $I=S^1$. Sometimes we will write $\length_X(c)$.
The space $X$ is geodesic if any two points $x,y\in X$ can be joined by a curve of length equal to $d(x,y)$.

For $s\geq 0$ the Hausdorff $s$-measure on a metric space $X$ is denoted by $\hm^s_X$ or simply by $\hm^s$. We choose the normalization constant in such a way that on Euclidean $\R^n$ the Hausdorff $n$-measure coincides with the Lebesgue measure. The Lebesgue measure of a subset $A\subset\R^n$ is denoted $|A|$. 

A map $\varphi\colon X\to Y$ between metric spaces $(X, d_X)$ and $(Y, d_Y)$ is called $(L, \alpha)$-H\"older continuous if $$d_Y(\varphi(x), \varphi(x'))\leq L\cdot d_X(x,x')^\alpha$$ for all $x,x'\in X$.

Let $L\geq 1$ and $\lambda_0>0$. A metric space $X$ is called $L$-Lipschitz $1$-connected up to scale $\lambda_0$ if every $\lambda$-Lipschitz curve $c\colon S^1\to X$ with $\lambda< \lambda_0$ extends to a $L\lambda$-Lipschitz map defined on all of $\overline{D}$. The space $X$ is called Lipschitz $1$-connected up to some scale if it is $L$-Lipschitz $1$-connected up to scale $\lambda_0$ for some $L\geq 1$ and $\lambda_0>0$.

\subsection{Ultralimits and asymptotic cones of metric spaces}\label{subsec:asymp-cones}

We briefly review some definitions and facts concerning ultralimits and asymptotic cones. For more details we refer for example to \cite{BrH99} or \cite{Dru00}. 

Let $\omega$ be a non-principal ultrafilter $\omega$ on $\N$, that is, a finitely additive measure on $\N$ such that every subset $A\subset\N$ is $\omega$-measurable with $\omega(A)\in \{0,1\}$ and such that $\omega(\N)=1$ and $\omega(A)=0$ whenever $A$ is finite. If $(Z,d)$ is a compact metric space then for every sequence $(z_n)\subset Z$ there exists a unique point $z\in Z$ such that 
$$\omega(\{n\in\N : d(z_n, z)>\varepsilon\}) = 0$$ for every $\varepsilon>0$. This point $z$ will be denoted by $\lim\nolimits_\omega z_n$ and we call it the ultralimit of the sequence $(z_n)$.

Let $(X_n, d_n, p_n)$ be a sequence of pointed metric spaces. We call a sequence of points $x_n\in X_n$ bounded if $$\sup_{n\in\N} d_n(x_n, p_n)<\infty.$$ Let us equip the set $\tilde{X}$ of bounded sequences in the sense above with the pseudo-metric $$\tilde{d}_\omega((x_n), (x'_n)):= \lim\nolimits_\omega d_n(x_n, x'_n).$$ The $\omega$-ultralimit of the sequence $(X_n, d_n, p_n)$ is the metric space obtained from $\tilde{X}$ by identifying points in $\tilde{X}$ of zero $d_\omega$-distance. We denote this space by $X_\omega = (X_n, d_n, p_n)_\omega$ and its metric by $d_\omega$. An element of $X_\omega$ will be denoted by $[(x_n)]$, where $(x_n)$ is an element of $\tilde{X}$. Ultralimits are always complete and ultralimits of sequences of geodesic metric spaces are again geodesic. 

Let $Y$ be a metric space, $\alpha\in(0,1]$ and $C\geq 0$. Let $\varphi_n\colon Y\to X_n$ be $(C,\alpha)$-H\"older maps, $n\in\N$. 
If the sequence $(\varphi_n)$ is bounded in the sense that
\begin{equation*}
 \sup_{n\in\N} d_n(\varphi_n(y), p_n)<\infty
\end{equation*}
for some and thus every $y\in Y$ then the assignment $y\mapsto [(\varphi_n(y))]$ defines a $(C,\alpha)$-H\"older map from $Y$ to $X_\omega$. We denote this map by $(\varphi_n)_\omega$ or $\lim\nolimits_\omega \varphi_n$. 

Let $(X,d)$ be a metric space, $(p_n)\subset X$ a sequence of basepoints and $(t_n)$ a sequence of positive real numbers satisfying $\lim_{n\to\infty} t_n=0$. The asymptotic cone of $X$ with respect to $(p_n)$, $(t_n)$ and $\omega$ is the $\omega$-ultralimit of the sequence $(X, t_n d, p_n)$. It will be denoted by $(X, t_n, p_n)_\omega$ or simply by $X_\omega$ if there is no danger of ambiguity. A geodesic metric space is Gromov hyperbolic if and only if every of its asymptotic cones is a metric tree, see \cite[Proposition 3.1.1]{Dru00}.

\subsection{The ${\rm CAT}(0)$ and ${\rm CAT}(0, \nu)$ conditions}\label{sec:coarse-CAT}

Let $(X,d)$ be a geodesic metric space. A geodesic triangle $\Delta$ in $X$ consists of three points in $X$ and a choice of three geodesics (the sides) connecting them. The sum of their lengths is called the perimeter of $\Delta$. Consider the comparison triangle $\overline{\Delta}\subset\R^2$ for $\Delta$. This is the unique (up to isometries) triangle in Euclidean $\R^2$ whose sides have the same lengths as the sides of $\Delta$. The triangle $\Delta$ is said to be ${\rm CAT}(0)$ if for all $x, y\in \Delta$ and their unique comparison points $\bar{x}, \bar{y}\in\overline{\Delta}$ the inequality $d(x,y) \leq |\bar{x} - \bar{y}|$ holds. In other words, $\Delta$ is at least as thin as $\overline{\Delta}$. If all geodesic triangles in $X$ are ${\rm CAT}(0)$ then $X$ is called ${\rm CAT}(0)$-space. We refer to \cite{BrH99} for details concerning the definitions above.

The following notion of coarse non-positive curvature was introduced by Gromov in \cite{Gro01}. See also \cite{DeG08}, \cite{Kar-Aditi-2011}, \cite{BuF12}, \cite{KS94} and the references therein.

\bd
 Let $(X,d)$ be a geodesic metric space and $\nu\geq 0$. A geodesic triangle $\Delta\subset X$ is called ${\rm CAT}(0,\nu)$ if for all $x,y\in \Delta$ the inequality $$d(x,y) \leq |\bar{x} - \bar{y}| + \nu$$ holds, where $\bar{x}, \bar{y}\in\overline{\Delta}$ are the comparison points in the comparison triangle $\overline{\Delta}\subset\R^2$.
\ed

We will need the following proposition whose simple proof is left to the reader. Compare with \cite[Proposition 3.2.8]{DeG08} and \cite[Theorem 8]{Kar-Aditi-2011}.

\bp\label{prop:ultralimit-CAT-consequence-sequence}
 Let $r_0>0$. Let $(X_n, d_n)$ be geodesic metric spaces, $n\in\N$, such that for every non-principal ultrafilter $\omega$ on $\N$ and every sequence of basepoints $p_n\in X_n$ all geodesic triangles of perimeter at most $r_0$ in the ultralimit $(X_n, d_n, p_n)_\omega$ are ${\rm CAT}(0)$. Then for every $\nu>0$ there exists $n_0\in\N$ such that if $n\geq n_0$ then all geodesic triangles in $X_n$ of perimeter at most $r_0$ are ${\rm CAT}(0,\nu)$.
\ep

\subsection{Sobolev maps with values in metric spaces}\label{sec:Sobolev-maps}

There exist several equivalent definitions of Sobolev maps from a Euclidean domain into a metric space, see e.g.~\cite{Amb90}, \cite{KS93}, \cite{Res97}, \cite{Res04}, \cite{Res06}, \cite{HKST01}, \cite{HKST15}, \cite{AT04}. We recall the definition from \cite{Res97} based on compositions with real-valued Lip\-schitz functions. We will only need Sobolev maps defined on the open unit disc $D$ of $\R^2$.

Let $(X,d)$ be a complete metric space and $p>1$. Let $L^p(D, X)$ be the set of measurable and essentially separably valued maps $u\colon D\to X$ such that for some and thus every $x\in X$ the function $u_x(z):= d(x, u(z))$ belongs to the classical space $L^p(D)$ of $p$-integrable functions on $D$. 

\bd
 A map $u\in L^p(D, X)$ belongs to the Sobolev space $W^{1,p}(D, X)$ if there exists $h\in L^p(D)$ such that $u_x$ is in the classical Sobolev space $W^{1,p}(D)$ for every $x\in X$ and its weak gradient satisfies $|\nabla u_x|\leq h$ almost everywhere. 
\ed

There are several natural notions of energy of a map $u\in W^{1,p}(D, X)$. Throughout this text we will use the Reshetnyak $p$-energy defined by $$E_+^p(u):= \inf\left\{\|h\|_{L^p(D)}^p\;\big|\; \text{$h$ as in the definition above}\right\}.$$

We will furthermore need the notion of trace of $u\in W^{1,p}(D, X)$. By \cite{KS93} there exists a representative $\bar{u}$ of $u$ such that the curve $t\mapsto \bar{u}(tv)$ with $t\in[1/2, 1)$ is absolutely continuous for almost every $v\in S^1$. The trace of $u$ is then defined by $$\trace(u)(v):= \lim_{t\nearrow 1}\bar{u}(tv)$$ for almost every $v\in S^1$. It can be shown that $\trace(u)\in L^p(S^1, X)$, see \cite{KS93}. If $u$ has a continuous extension $\hat{u}$ to $\overline{D}$ then $\trace(u)$ is just the restriction of $\hat{u}$ to ${S^1}$. 

Every map $u\in W^{1,p}(D, X)$ has an approximate metric derivative at almost every point $z\in D$ in the following sense, see \cite{Kar07} and \cite{LW15-Plateau}. There exists a unique seminorm on $\R^2$, denoted $\apmd u_z$, such that 
 \begin{equation*}
    \ap\lim_{z'\to z}\frac{d(u(z'), u(z)) - \apmd u_z(z'-z)}{|z'-z|} = 0,
 \end{equation*}
where $\ap\lim$ denotes the approximate limit, see \cite{EG92}. If $u$ is Lipschitz then the approximate limit can be replaced by an honest limit. It follows from \cite{LW15-Plateau} that $$E_+^p(u) = \int_D\mathcal{I}_+^p(\apmd u_z)\,dz,$$ where for a seminorm $s$ on $\R^2$ we have set $\mathcal{I}_+^p(s):= \max\{s(v)^p : |v|=1\}$.

\bd
The (parameterized Hausdorff) area of a map $u\in W^{1,2}(D, X)$ is defined by $$\Area(u):= \int_D \jac(\apmd u_z)\,dz,$$ where the Jacobian $\jac(s)$ of a seminorm $s$ on $\R^2$ is the Hausdorff $2$-measure in $(\R^2, s)$ of the Euclidean unit square if $s$ is a norm and $\jac(s)=0$ otherwise. 
\ed

If $u\in W^{1,2}(D, X)$ satisfies Lusin's property (N), thus sends sets of Lebesgue measure zero to sets of Hausdorff $2$-measure zero, then $$\Area(u)=  \int_X\#\{z : u(z) = x\} \,d\hm^2(x)$$ by the area formula \cite{Kir94}. In particular, if $u$ is injective then $\Area(u) = \hm^2(u(D))$. The area and energy are related by $\Area(u)\leq E_+^2(u)$ for every $u\in W^{1,2}(D,X)$, see \cite[Lemma 7.2]{LW15-Plateau}.

We will need the following simple observation.

 \bl\label{lem:product-Sobolev-maps}
  Let $X$ and $Y$ be complete metric spaces and let $u\in W^{1,2}(D, X)$ and $v\in W^{1,2}(D, Y)$. Then the map $w = (u,v)$ belongs to $W^{1,2}(D, X\times Y)$ and its approximate metric derivative at almost every $z\in D$ satisfies $$\left[\apmd w_z(\xi)\right]^2 = \left[\apmd u_z(\xi)\right]^2 + \left[\apmd v_z(\xi)\right]^2$$ for all $\xi\in\R^2$.
 \el
 
 \begin{proof}
  It follows for example from \cite[Proposition 3.2]{LW15-Plateau} that $w \in W^{1,2}(D, X\times Y)$. Denote the distances on $X$, $Y$, and $X\times Y$ by $d_X$, $d_Y$, and $d_{X\times Y}$, respectively. Let $z\in D$ be a point at which each of the maps $u$, $v$, and $w$ is approximately metrically differentiable and let $\xi\in S^1$. Thus, there exist a sequence $(r_n)$ of positive real numbers and a sequence $(\xi_n)\subset \R^2$ with $r_n\to 0$ and $\xi_n\to \xi$ and such that $$\lim_{n\to\infty}r_n^{-1}\cdot d_X(u(z + r_n\xi_n), u(z)) = \apmd u_z(\xi),$$
  $$\lim_{n\to\infty}r_n^{-1}\cdot d_Y(v(z + r_n\xi_n), v(z)) = \apmd v_z(\xi),$$ $$\lim_{n\to\infty}r_n^{-1}\cdot d_{X\times Y}(w(z + r_n\xi_n), w(z)) = \apmd w_z(\xi),$$ from which we infer that
 \begin{equation*}
 \begin{split} 
  [\apmd w_z(\xi)]^2 &= \lim_{n\to\infty} r_n^{-2}\cdot \left[d_X(u(z + r_n\xi_n), u(z))^2 +  d_Y(v(z + r_n\xi_n), v(z))^2 \right]\\
  & = \left[\apmd u_z(\xi)\right]^2 + \left[\apmd v_z(\xi)\right]^2.
 \end{split}
 \end{equation*}
 This completes the proof.
 \end{proof}

\subsection{Sobolev Dehn function}\label{sec:Sobolev-Dehn}

We will need the following variant of the Lipschitz Dehn function introduced in \cite{LWY16}. The (Sobolev) filling area of a Lipschitz curve $c\colon S^1\to X$ in a complete metric space $X$ is defined by
\begin{equation*}
 \Fillarea(c):= \inf\left\{\Area(u): u\in W^{1,2}(D, X), \trace(u) = c\right\}.
\end{equation*}
Sometimes we will write $\Fillarea_X(c)$ and similarly for the Lipschitz filling area.
The (Sobolev) Dehn function is given by
\begin{equation*}
 \delta_X(r):= \sup\left\{\Fillarea(c) : \text{$c\colon S^1\to X$ is Lipschitz, $\length(c)\leq r$}\right\}
\end{equation*}
for all $r\geq 0$.
We clearly have $\delta_X(r)\leq \deltalip_X(r)$ for all $r$. Moreover, equality holds for example if $X$ is geodesic and Lipschitz $1$-connected up to some scale, see \cite[Proposition 3.1]{LWY16}. 

One of the principal advantages the Sobolev Dehn function has over its more classical Lipschitz analog is the following stability property whose proof is the same as that of \cite[Corollary 5.3]{LWY16}. 

\bt\label{thm:Sobolev-Dehn-stable}
 Let $0< r_0\leq\infty$ and $C>0$. Let $(\varepsilon_n)$ be a sequence of non-negative real numbers tending to $0$. For every $n\in\N$ let $(X_n, d_n)$ be a proper, geodesic metric space satisfying 
 \begin{equation}\label{eq:quad-isop-error}
  \delta_{X_n}(r)\leq (C + \varepsilon_n) \cdot r^2 + \varepsilon_n
 \end{equation}
  for all $r\in(0,r_0)$. Then for every non-principal ultrafilter $\omega$ on $\N$ and every sequence of basepoints $p_n\in X_n$ the ultralimit $X_\omega=(X_n, d_n, p_n)_\omega$ satisfies $$\delta_{X_\omega}(r)\leq C\cdot r^2$$ for all $r\in(0,r_0)$.
\et

We mention that rescalings of spaces with asymptotically quadratic Dehn function satisfy the bound \eqref{eq:quad-isop-error}. More precisely, let $C>0$ and let $(X, d)$ be a metric space such that $$\limsup_{r\to\infty}\frac{\delta_X(r)}{r^2}\leq C.$$ Let $(t_n)$ be a sequence of positive real numbers tending to $0$ and let $X_n$ be the metric space given by $X_n = (X, t_n  d)$. Then there exists a sequence $(\varepsilon_n)$ of positive real numbers tending to $0$ such that $$\delta_{X_n}(r)\leq (C+\varepsilon_n)\cdot r^2 + \varepsilon_n$$ for all $r\geq 0$ and all $n\in\N$.

We end the section with the following result which guarantees the existence of suitable thickenings of a metric space.

\bp\label{prop:good-thickening}
 There exists a universal constant $L\geq 1$ with the following property. Let $0<r_0\leq \infty$ and $C >0$. Let $\varepsilon\in(0,1)$ and let $X$ be a proper, geodesic metric space satisfying $$\delta_X(r) \leq C\cdot r^2 + \varepsilon^2$$ for all $r\in(0,r_0)$. Then there exists a proper, geodesic metric space $Y$ with the following properties:
  \begin{enumerate}
   \item $Y$ contains $X$ and lies at Hausdorff distance at most $\varepsilon$ from $X$.
   \item $Y$ is $L$-Lipschitz $1$-connected up to scale $L^{-1}\varepsilon$.
   \item $Y$ satisfies $\delta_Y(r)\leq (C+L^2)\cdot r^2$ for all $r\in(0,r_0)$ and $$\delta_Y(r) \leq \left(C + \sqrt{\varepsilon}\right) \cdot r^2$$ for all $r\in(0,r_0)$ with $r\geq L\sqrt{\varepsilon}$.
  \end{enumerate}
\ep

\begin{proof}
  This follows as in the proof of \cite[Proposition 3.5]{LWY16}.
\end{proof}

\section{Spaces with property (ET)}\label{sec:ET}

Recall from \cite{LW15-Plateau} that a complete metric space $X$ is said to have property (ET) if for every $u\in W^{1,2}(D, X)$ the approximate metric derivative $\apmd u_z$ comes from a possibly degenerate inner product at almost every $z\in D$.

The aim of this section is to establish the following result, which generalizes \cite[Theorem 5.2]{LW-isoperimetric} from the setting of proper to that of complete metric spaces.
 
 \bt\label{thm:eucl-isop-property-ET}
  Let $X$ be a complete, geodesic metric space and $r_0>0$. If $\delta_X(r)\leq \frac{r^2}{4\pi}$ for all $r\in(0,r_0)$ then $X$ has property (ET).
 \et
 
The proof of this theorem is similar to the proof of \cite[Theorem 5.1]{Wen08-sharp}. We first show:
 
 \bl\label{lem:non-ET}
  Let $(X, d)$ be a complete metric space. If $X$ does not have property (ET) then there exists a non-Euclidean norm $\|\cdot\|$ on $\R^2$ with the following properties. For every finite set $\{v_1,\dots, v_n\}\subset\R^2$ and every $\lambda>1$ there exist $\delta>0$ arbitrarily small and points $x_1,\dots, x_n\in X$ such that 
 \begin{equation}\label{eq:bilip-ineq-X-V}
 \lambda^{-1}\delta\cdot\|v_k - v_m\|\leq d(x_k, x_m) \leq \lambda \delta\cdot\|v_k-v_m\|
 \end{equation}
  for all $k,m=1,\dots, n$.
 \el

\begin{proof}
 Since $X$ does not have property (ET) there exists $u\in W^{1,2}(D, X)$ whose approximate metric derivative $\apmd u_z$ does not come from a possibly degenerate inner product almost everywhere. 
By \cite[Proposition 4.3]{LW15-Plateau} there thus exists a measurable subset $K\subset D$ of strictly positive measure with the following properties. Firstly, the approximate metric derivative $\apmd u_z$ exists for every $z\in K$ and is a non-Euclidean norm. Secondly, for every $z\in K$ and $\lambda>1$ the norm $\apmd u_z$ satisfies $$\lambda^{-1}\cdot \apmd u_z(z' - z'') \leq d(u(z'), u(z'')) \leq \lambda\cdot \apmd u_z(z'-z'')$$ for all $z', z''\in K$ contained in a sufficiently small ball around $z$. 

Fix a Lebesgue density point $z\in K$ of $K$ and set $\|\cdot\|:= \apmd u_z$. We may assume that $z=0$. Let $\{v_1,\dots, v_n\}\subset \R^2$ be a finite set and $\lambda>1$. Then any ball around $v_k$ intersects the set $\delta^{-1} K= \{\delta^{-1} x: x\in K\}$ for every sufficiently small $\delta>0$, depending on the radius of the ball. If $z_k\in \delta^{-1} K$ is sufficiently close to $v_k$ then the points $x_k = u(\delta z_k)$ satisfy \eqref{eq:bilip-ineq-X-V}. This concludes the proof. 
\end{proof}

 \begin{proof}[Proof of Theorem~\ref{thm:eucl-isop-property-ET}]
 We denote the metric on $X$ by $d$. We argue by contradiction and assume that $X$ does not have property (ET). Let $\|\cdot \|$ be a non-Euclidean norm given by Lemma~\ref{lem:non-ET} and denote by $V$ the normed space $(\R^2, \|\cdot\|)$. Let $\mathbb{I}_V\subset V$ be an isoperimetric set for $V$. Thus, $\mathbb{I}_V$ is a convex subset of largest area among all convex subsets of $V$ with given boundary length. Since $V$ is not Euclidean we have
 \begin{equation}\label{eq:isop-non-Euclidean}
 \hm_V^2(\mathbb{I}_V) > \frac{1}{4\pi}\cdot \length_V(\partial \mathbb{I}_V)^2,
\end{equation}
  see for example \cite[Lemma 5.1]{LW-isoperimetric}. 
 
Let $\gamma\colon S^1\to V$ be a constant speed parametrization of $\partial \mathbb{I}_V$. Let $\lambda>1$ be sufficiently close to $1$ and let $n\in \N$ be sufficiently large, to be determined later. For $k=1,\dots, n$ define $z_k:= e^{2\pi i \frac{k}{n}}$ and $v_k:= \gamma(z_k)$. By Lemma~\ref{lem:non-ET} there exist $\delta>0$ arbitrarily small and $x_1,\dots, x_n\in X$ such that $$\lambda^{-1}\delta \cdot\|v_k - v_m\|\leq d(x_k, x_m) \leq \lambda \delta\cdot\|v_k-v_m\|$$ for all $k,m$. After replacing the norm $\|\cdot\|$ by the rescaled norm $\delta\|\cdot\|$ we may assume that $\delta = 1$ and that $\lambda\cdot \length_V(\partial \mathbb{I}_V)< r_0$.

Let $c\colon S^1\to X$ be the curve satisfying $c(z_k) = x_k$ and which is geodesic on the segment of $S^1$ between $z_k$ and $z_{k+1}$. Notice that $\length(c) \leq \lambda\cdot \length_V(\partial \mathbb{I}_V)$. By the Euclidean isoperimetric inequality on $X$ up to scale $r_0$ there exists $u\in W^{1,2}(D, X)$ such that $\trace(u) = c$ and $$\Area(u)\leq \frac{1}{4\pi}\cdot \length(c)^2 \leq \frac{\lambda^2}{4\pi}\cdot \length_V(\partial \mathbb{I}_V)^2.$$
View $V$ as a linear subspace of the space $\ell^\infty$ of bounded sequences in $\R$ with the supremum norm. Since $\ell^\infty$ is an injective metric space there exists a $\lambda$-Lipschitz map $\varphi\colon X\to \ell^\infty$ extending the map which sends  $x_k$ to $v_k$ for every $k$.
Then the map $\varphi\circ u$ belongs to $W^{1,2}(D, \ell^\infty)$ and satisfies $\trace(\varphi\circ u) = \varphi\circ c$ as well as $$\Area(\varphi\circ u) \leq \frac{\lambda^4}{4\pi}\cdot \length_V(\partial \mathbb{I}_V)^2.$$ By \cite[Proposition 3.1]{LWY16} there exists for every $\varepsilon>0$ a Lipschitz map $v\colon\overline{D}\to \ell^\infty$ with $v|_{S^1} = \varphi\circ c$ and $\Area(v)\leq \Area(\varphi\circ u) + \varepsilon$. 

We can connect the curves $\varphi\circ c$ and $\gamma$ by a Lipschitz homotopy $\varrho\colon S^1\times[0,1]\to \ell^\infty$ of small area as follows. Let $\varrho(z,0) = \varphi(c(z))$ and $\varrho(z,1) = \gamma(z)$ and let $\varrho(z_k,t) = v_k$ for all $k$ and every $t\in[0,1]$. The restriction of $\varrho$ to the boundary $\partial A_k$ of 
$$A_k:= \left\{e^{2\pi i \frac{\theta}{n}}: \theta\in [k, k+1]\right\}\times [0,1]$$ is a Lipschitz curve of length $$\length(\varrho|_{\partial A_k})\leq \frac{1+ \lambda^2}{n}\cdot \length_V(\partial \mathbb{I}_V).$$ By the quadratic isoperimetric inequality in $\ell^\infty$, there thus exists a Lipschitz extension $\varrho$ of $\varrho|_{\partial A_k}$ to $A_k$ with area  $$\Area(\varrho|_{A_k})\leq \frac{C(1+ \lambda^2)^2}{n^2}\cdot \length_V(\partial\mathbb{I}_V)^2,$$ where $C$ is a constant. This defines $\varrho$ on all of $S^1\times[0,1]$ and thus provides a Lipschitz homotopy from $\varphi\circ c$ to $\gamma$ satisfying
\begin{equation*}
  \Area(\varrho) \leq \frac{C(1+\lambda^2)^2}{n}\cdot \length_V(\partial \mathbb{I}_V)^2.
\end{equation*}
Finally, we can construct a Lipschitz map $w\colon \overline{D}\to \ell^\infty$ with $w|_{S^1}=\gamma$ and 
\begin{equation}\label{eq:area-w} 
\Area(w)\leq \left[\frac{\lambda^4}{4\pi} + \frac{C(1+\lambda^2)^2}{n}\right]\cdot \length_V(\partial \mathbb{I}_V)^2 + \varepsilon
\end{equation}
by gluing $v$ and $\varrho$ along $S^1$ and $S^1\times\{0\}$.
By the quasi-convexity of the Hausdorff $2$-measure \cite{BI12} we have $\hm_V^2(\mathbb{I}_V)\leq \Area(w)$. For $\varepsilon>0$ sufficiently small, $\lambda>1$ sufficiently close to $1$, and $n$ sufficiently large, inequality \eqref{eq:area-w} thus contradicts \eqref{eq:isop-non-Euclidean}. This completes the proof.
\end{proof}

We will furthermore need the following proposition which is a direct consequence of Lemma~\ref{lem:product-Sobolev-maps}.

\bp\label{prop:property-ET-products}
 Let $X$ and $Y$ be complete metric spaces. If $X$ and $Y$ have property (ET) then the space $X\times Y$, equipped with the Euclidean product metric, also has property (ET).
\ep

 \section{Discretization of fillings in spaces with property (ET)}\label{sec:discretization}
 
The following theorem will be one of the main ingredients in the proof of the existence of energy and area minimizers in ultralimits and asymptotic cones. It will be used in Section~\ref{sec:competitors} to bound the filling area of approximating curves by the filling area of the limit curve in an ultralimit.

In what follows, two disjoint Jordan curves $\Gamma, \Gamma'\subset\R^2$ are called separated if also their Jordan domains are disjoint.
 
 \bt\label{thm:discretization}
  Let $X$ be a complete, geodesic metric space with property (ET) and let $c\colon S^1\to X$ be a Lipschitz curve with $\Fillarea(c)<\infty$. Then for every $\varepsilon>0$ and every sufficiently large $\lambda\geq 1$ there exist a finite collection $\{\Gamma_1,\dots, \Gamma_k\}$ of pairwise separated smooth convex Jordan curves $\Gamma_i\subset D$ and a $\lambda$-Lipschitz extension $\varphi\colon K\to X$ of $c$ to the set $K:= S^1\cup\Gamma_1\cup\dots\cup \Gamma_k$ such that $$\sum_{i=1}^k \length\left(\varphi|_{\Gamma_i}\right)^2 \leq 4\pi\cdot \Fillarea(c) + \varepsilon$$ and 
 \begin{equation}\label{eq:thm-discreti-small-rest} 
 \left| D\setminus \bigcup_{i=1}^k\Omega_i\right| \leq \frac{\varepsilon}{\lambda^2},
\end{equation}
 where $\Omega_i$ denotes the Jordan domain enclosed by $\Gamma_i$. 
 \et

 The rest of this section is devoted to the proof of the theorem above. From now on, the metric on $X$ will be denoted by $d$. We need the following simple observation.

 \bl\label{lem:ext-smooth-Jordan-Lip}
  For every smooth Jordan curve $\Gamma\subset \R^2$ there exists $\nu>0$ with the following property. If $A\subset\Gamma$ is a finite $\nu$-dense set in $\Gamma$ and $\varphi\colon A\to X$ a $\lambda$-Lipschitz map to a geodesic metric space $X$ then the piecewise geodesic extension $\bar{\varphi}\colon \Gamma\to X$ of $\varphi$ is $3\lambda$-Lipschitz.
 \el
 
 Here, $A$ and $\Gamma$ are equipped with the Euclidean metric from $\R^2$.
 
\begin{proof}
 Let $d_\Gamma$ denote the length metric on $\Gamma$. Notice that the identity map from $(\Gamma, d_\Gamma)$ to $(\Gamma, |\cdot |)$ is $1$-Lipschitz. Moreover, since $\Gamma$ is smooth there exists $\nu>0$ such that $d_\Gamma(x,y) \leq 3\cdot |x-y|$ for all $x,y\in\Gamma$ with $|x-y|\leq 4\nu$.
 
 Now, let $A\subset\Gamma$ be a finite $\nu$-dense subset of $\Gamma$ and let $\varphi\colon A\to X$ be a $\lambda$-Lipschitz map. Denote by $\bar{\varphi}\colon \Gamma\to X$ the piecewise geodesic extension. Then $\bar{\varphi}$ is $\lambda$-Lipschitz as a map from $(\Gamma, d_\Gamma)$ to $X$. Let $x,y\in\Gamma$. If $|x-y|\leq 4\nu$ then $$d(\bar{\varphi}(x), \bar{\varphi}(y))\leq \lambda d_\Gamma(x,y) \leq 3\lambda|x-y|.$$ If $|x-y|\geq 4\nu$ then let $a,b\in A$ be nearest points for $x$ and $y$, respectively. Then $|x-a| + |y-b|\leq 2\nu \leq \frac{1}{2} \cdot |x-y|$ and thus $|a-b| \leq \frac{3}{2}\cdot |x-y|$. Hence
\begin{equation*}
 \begin{split}
  d(\bar{\varphi}(x), \bar{\varphi}(y)) &\leq d(\bar{\varphi}(x), \bar{\varphi}(a)) + d(\varphi(a), \varphi(b)) + d(\bar{\varphi}(b), \bar{\varphi}(y))\\
  & \leq 3\lambda  |x-a| + \lambda |a-b| + 3\lambda |b-y|\\
  & \leq 3\lambda  |x-y|.
 \end{split}
\end{equation*}
This completes the proof.
 \end{proof}

Let $X$ and $c$ be as in the statement of Theorem~\ref{thm:discretization} and let $\varepsilon>0$. Let $u\in W^{1,2}(D, X)$ be such that $\trace(u)= c$ and $\Area(u) \leq \Fillarea(c) + \varepsilon$. We will use the map $u$ to construct curves $\Gamma_i$ and a $\lambda$-Lipschitz map $\varphi\colon S^1\cup\Gamma_1\cup\dots\cup \Gamma_k\to X$ which extends $c$ and satisfies \eqref{eq:thm-discreti-small-rest} and $$\sum_{i=1}^k \length\left(\varphi|_{\Gamma_i}\right)^2 \leq 4\pi\cdot \Area(u) + \varepsilon.$$
Roughly speaking, the $\Gamma_i$ will be suitably chosen ellipses and $\varphi$ will coincide with $u$ on a finite and sufficiently dense set of points in each $\Gamma_i$ and will be piecewise geodesic on each $\Gamma_i$. We first show that we may assume $u$ to have some additional properties. These will be used in the proofs of Propositions~\ref{prop:Lip-big-set} and \ref{prop:almost-isom}.
 
 \bl\label{lem:add-props}
  We may assume that the restriction of $u$ to $D\setminus B(0,\frac{1}{2})$ is Lipschitz continuous and that there exists $\delta>0$ such that $$d(u(z), u(z'))\geq  \delta \cdot |z-z'|$$ for all $z,z'\in D$. 
 \el

 \begin{proof}
  Let $u$ and $c$ be as in the paragraph preceding the lemma. Define a map $u'\colon D\to X$ by $u'(z)= u(2z)$ if $|z|<\frac{1}{2}$ and $u'(z) = c(z/|z|)$ if $|z|\geq \frac{1}{2}$. Then $u'$ belongs to $W^{1,2}(D, X)$ by \cite[Theorem 1.12.3]{KS93}. Moreover, $u'$ is Lipschitz on $D\setminus B(0,\frac{1}{2})$ and satisfies $\trace(u') = c$ and $\Area(u') = \Area(u)$.

Now, let $Y$ be the space $X\times \R^2$ equipped with the Euclidean product metric, which we denote by $d_Y$. Then $Y$ is complete, geodesic, and has property (ET) by Proposition~\ref{prop:property-ET-products}. For $\delta>0$ consider the map $u_\delta\colon D\to Y$ given by $u_\delta(z):= (u'(z), \delta z)$. Then $u_\delta$ belongs to $W^{1,2}(D, Y)$, is Lipschitz on $D\setminus B(0,\frac{1}{2})$ and satisfies $$d_Y(u_\delta(z), u_\delta(z'))\geq \delta\cdot |z-z'|$$ for all $z,z\in D$. Moreover, the trace of $u_\delta$ is the curve given by $c_\delta(z)=(c(z), \delta z)$ for all $z\in S^1$. 
We claim that $\Area(u_\delta) \to \Area(u)$ as $\delta\to 0$. For this, first note that $$(\apmd (u_\delta)_z(w))^2 = (\apmd u_z(w))^2 + \delta^2\cdot|w|^2$$ for almost every $z\in D$ and every $w\in \R^2$ by Lemma~\ref{lem:product-Sobolev-maps}. Hence, $\jac(\apmd (u_\delta)_z)$ converges to $\jac(\apmd u_z)$ and $$\jac(\apmd u_z)\leq \jac(\apmd (u_\delta)_z) \leq \mathcal{I}_+^2(\apmd u_z) + \delta^2$$ for almost every $z\in D$.  The dominated convergence theorem now implies the claim. 

Choose $\delta>0$ so small that $\Area(u_\delta)<\Area(u) + \varepsilon$. Let $P\colon Y\to X$ be the natural projection and note that $P$ is $1$-Lipschitz. Suppose that we can use the map $u_\delta$ to construct, for every $\lambda\geq 1$ sufficiently large, a finite collection $\{\Gamma_1,\dots, \Gamma_k\}$ of pairwise separated smooth convex Jordan curves $\Gamma_i\subset D$ and a $\lambda$-Lipschitz extension $\varphi'\colon K\to Y$ of $c_\delta$ to the set $K:= S^1\cup\Gamma_1\cup\dots\cup \Gamma_k$ such that \eqref{eq:thm-discreti-small-rest} holds and $$\sum_{i=1}^k \length\left(\varphi'|_{\Gamma_i}\right)^2 \leq 4\pi\cdot \Area(u_\delta) + \varepsilon.$$ Then the map $\varphi:= P\circ \varphi'$ is a $\lambda$-Lipschitz extension of $c$ and satisfies 
\begin{equation*}
 \sum_{i=1}^k \length\left(\varphi|_{\Gamma_i}\right)^2 \leq 4\pi\cdot \Area(u_\delta) + \varepsilon \leq 4\pi\cdot \Fillarea(c) + (8\pi + 1)\varepsilon.
 \end{equation*}
This shows that it is indeed enough to use the map $u_\delta$ to construct the desired Jordan curves $\Gamma_i$ and the Lipschitz map $\varphi$.
Since $u_\delta$ has all the properties in the statement of the lemma and since $Y$ has the same properties as $X$ the proof is complete.
 \end{proof}
 
 From now on, we assume that $u$ also satisfies the properties of Lemma~\ref{lem:add-props}. 
 
 \bp\label{prop:Lip-big-set}
  For every sufficiently large $\lambda\geq 1$ there exists a measurable set $E\subset B(0,\frac{3}{4})$ with $|E|<\frac{\varepsilon}{\lambda^2}$ and such that the restriction of $u$ to $D\setminus E$ is $\lambda$-Lipschitz.
 \ep
 
 This essentially follows from the proof of \cite[Theorem 8.2.1]{HKST15}. For the convenience of the reader, we provide the proof.

 \begin{proof}
  Define a function by $h(z):= \mathcal{I}_+^1(\apmd u_z)$ and note that $h\in L^2(D)$ and that $h\leq L$ on $D\setminus \bar{B}(0,\frac{1}{2})$ for some $L>0$ since $u$ is Lipschitz on this set. By the proof of \cite[Proposition 3.2]{LW15-Plateau} there exists a set $N\subset D$ of measure zero and a constant $C>0$ such that $$d(u(z), u(z'))\leq |z-z'|\cdot (g(z) + g(z'))$$ for all $z,z'\in D\setminus N$, where $g(z):= C\cdot M(h)(z)$ and $M(h)$ denotes the maximal function of $h$. Notice that $g\in L^2(D)$ by the maximal function theorem \cite[Theorem 2.2]{Hei01}.
  
For $\lambda>0$ set $$E_\lambda:= \left\{z\in D: g(z) >\frac{\lambda}{2}\right\}$$ and observe that $u$ is $\lambda$-Lipschitz on $D\setminus (E_\lambda \cup N)$. Chebyshev's inequality implies $$|E_\lambda| \leq \frac{4}{\lambda^2}\cdot \int_{E_\lambda}g^2(z)\,dz$$ and hence, by the absolute continuity of the integral, we have $|E_\lambda|< \frac{\varepsilon}{\lambda^2}$ whenever $\lambda\geq 1$ is sufficiently large. 

Finally, since $h\leq L$ on $D\setminus \bar{B}(0,\frac{1}{2})$ a direct calculation shows that $g$ is bounded on $D\setminus B(0,\frac{3}{4})$ by a constant depending only on $C$, $L$ and $\|h\|_{L^1(D)}$. In particular, for every sufficiently large $\lambda\geq 1$ we have $E_\lambda\subset B(0,\frac{3}{4})$. The continuity of $u$ on $D\setminus \bar{B}(0, \frac{1}{2})$ now implies that $u$ is $\lambda$-Lipschitz on the set $D\setminus E$, where $E= E_\lambda \cup (N\cap B(0, \frac{3}{4}))$. This completes the proof.
 \end{proof}
 
Let $\lambda\geq 1$ be sufficiently large and set $F:= D\setminus E$, where $E$ is as in the proposition above. By the proposition and lemma above, the restriction of $u$ to $F$ is biLipschitz. Since $F$ contains the annulus $D\setminus B(0,\frac{3}{4})$, the map $u$ extends to a $\lambda$-Lipschitz map on $F\cup S^1$. We denote the extension by $u$ again and notice that, by the definition of trace, we have $u|_{S^1}= c$.

 \bp\label{prop:almost-isom}
  There exist pairwise disjoint compact subsets $K_1,\dots, K_m\subset F$ and inner product norms $\|\cdot\|_i$, $i=1,\dots, m$, such that $|F\setminus \cup_{i=1}^m K_i|<\frac{\varepsilon}{\lambda^2}$ and $$(1+\varepsilon)^{-1}\cdot \|z-z'\|_i \leq d(u(z), u(z'))\leq (1+\varepsilon)\cdot \|z-z'\|_i$$ for all $z,z'\in K_i$ and for all $i=1,\dots, m$.
 \ep
 
From the area formula we conclude that $$\hm^2_{Y_i} (K_i) \leq (1+\varepsilon)^{2}\cdot  \Area(u|_{K_i}),$$ where we have set $Y_i:= (\R^2, \|\cdot\|_i)$.

 \begin{proof}
  Since $X$ has property (ET) and the restriction of $u$ to $F$ is $L$-biLipschitz for some $L\geq 1$ it follows that the approximate metric derivative $\apmd u_z$ comes from a (non-degenerate) inner product for almost every $z\in F$ and  satisfies $$L^{-1} \cdot |v| \leq \apmd u_z(v) \leq L\cdot |v|$$ for every $v\in \R^2$. Now, the proposition follows for example from \cite[Proposition 4.3]{LW15-Plateau} and the inner regularity of the Lebesgue measure.
 \end{proof}
 
 The next proposition is essentially a consequence of the Vitali covering theorem and will be applied to the sets $K_i$ from above, viewed as subsets of $(\R^2, \|\cdot \|_i)$.
 
 \bp\label{prop:covering}
  Let $\|\cdot\|$ be an inner product norm on $\R^2$, let $K$ be a bounded and measurable subset of $Y = (\R^2, \|\cdot\|)$ of positive Lebesgue measure and $U\subset Y$ open with $\overline{K}\subset U$. Then for all $\varepsilon', \rho>0$ there exists a finite collection of closed balls $\bar{B}_j= \bar{B}_Y(x_j, r_j)\subset U$, $j=1,\dots, k$, with the following properties:
  \begin{enumerate}
   \item The balls $\bar{B}_Y(x_j, (1+\varepsilon')r_j)$ are pairwise disjoint.
   \item For every $j$ the set $K\cap \partial \bar{B}_j$ is $\rho r_j$-dense in $\partial \bar{B}_j$.
     \item $\hm^2_{Y}\left(K\setminus \cup_{j=1}^{k} \bar{B}_j\right) \leq 6\varepsilon'\cdot \hm^2_{Y}(K)$.
  \end{enumerate}
 \ep
 
 \begin{proof}
  We may assume that $\|\cdot\|$ is the standard Euclidean norm $|\cdot |$ and that $\varepsilon'<\frac{1}{6}$ and $\rho<1$. Let $K'\subset K$ be the set of Lebesgue density points of $K$ and notice that $|K\setminus K'| = 0$. Let $V\subset U$ be a suitably chosen open bounded set containing $K'$ and $0<\mu<1$ sufficiently close to $1$, both to be determined later. Let $\mathfrak{B}$ be the family of all closed balls $\bar{B}(x,s)\subset V$ with $x\in K'$ and $s>0$ and such that $$|K'\cap \bar{B}(x,s)| \geq \mu \pi s^2.$$ Then $\mathfrak{B}$ is a fine covering of $K'$ in the sense of Vitali. Thus, by the Vitali covering theorem \cite[Theorem 2.8]{Mat95}, there exists a finite collection of closed pairwise disjoint balls $\bar{B}(x_j, s_j)\in \mathfrak{B}$, $j=1,\dots, k$, such that $$\left|K'\setminus \bigcup_{j=1}^k \bar{B}(x_j, s_j)\right|< \frac{\varepsilon'}{2} \cdot |V|.$$
  
Set $t:= (1+2\varepsilon')^{-1}$ and $t':=  (1+\varepsilon')^{-1}$. We claim that there exists $ts_j<r_j<t's_j$ such that $K'\cap \partial\bar{B}(x_j, r_j)$ is $\rho r_j$-dense in $\partial\bar{B}(x_j, r_j)$. Suppose this is not true. Then $$\hm^1\left(K'\cap \partial\bar{B}(x_j, r)\right) \leq (2\pi - 2\rho)r = 2\pi(1-\rho') r$$ for every $r\in(ts_j, t's_j)$, where $\rho' = \frac{\rho}{\pi}$. Hence, we obtain
\begin{equation*}
 \begin{split}
  \left|K'\cap \bar{B}(x_j, s_j)\right| &\leq \left|\bar{B}(x_j, ts_j)\right| + \left|\bar{B}(x_j, s_j)\setminus \bar{B}(x_j, t's_j)\right| + \int_{ts_j}^{t's_j} 2\pi(1-\rho') r\,dr\\
  &=  \pi s_j^2\cdot\left[1-\rho'(t'^2 - t^2)\right].
 \end{split}
\end{equation*}
However, if $\mu$ had been chosen sufficiently close to $1$ only depending on $\varepsilon'$ and $\rho$ then this is strictly smaller than $\mu\pi s_j^2$, which is a contradiction. This proves our claim.
Since $(1+\varepsilon') r_j< s_j$ it follows that the balls $\bar{B}(x_j, (1+\varepsilon')r_j)$ are pairwise disjoint and that $K\cap \partial \bar{B}(x_j, r_j)$ is $\rho r_j$-dense in $\partial\bar{B}(x_j, r_j)$.

 It remains to prove that property (iii) of the proposition holds. Notice that $$\left| \bar{B}(x_j, s_j)\setminus \bar{B}(x_j, ts_j)\right| = \pi s_j^2\cdot(1- t^2)< 5\varepsilon' \cdot |\bar{B}(x_j, s_j)|$$ and hence 
\begin{equation*}
 \begin{split}
  \left| K\setminus \bigcup_{j=1}^k \bar{B}(x_j, r_j)\right| &\leq \left|K'\setminus \bigcup_{j=1}^k \bar{B}(x_j, s_j)\right| + \sum_{j=1}^k \left| \bar{B}(x_j, s_j)\setminus \bar{B}(x_j, ts_j)\right|\\
 &< \frac{\varepsilon'}{2}\cdot|V| + 5\varepsilon' \cdot \left|\bigcup_{j=1}^k \bar{B}(x_j, s_j)\right|\\
 &\leq \frac{11}{2}\varepsilon'\cdot |V|.
 \end{split}
\end{equation*}
Since $|V|$ can be chosen arbitrarily close to $|K|$ property (iii) follows. This completes the proof.
 \end{proof}
 
 Let $K_i$ and $\|\cdot\|_i$ be as in Proposition~\ref{prop:almost-isom}. We may assume that each $K_i$ has positive measure. Let $\nu>0$ be so small that Lemma~\ref{lem:ext-smooth-Jordan-Lip} applies with $\nu$ to the boundary of the unit ball of $Y_i$, viewed as a subset of $(\R^2, |\cdot|)$, for every $i$. Recall that $Y_i = (\R^2, \|\cdot\|_i)$. Let $L\geq 1$ be such that the identity map from $(\R^2, |\cdot|)$ to $Y_i$ is $L$-biLipschitz for all $i$. 
 Choose open sets $U_i\subset\R^2$ satisfying $K_i \subset U_i\subset \overline{U}_i\subset D$ and $$\hm^2_{Y_i}(U_i) \leq (1+\varepsilon)\cdot \hm^2_{Y_i}(K_i)$$ and such that $\dist(U_i, U_j)>0$ for all $i\not=j$. Let $\rho_0>0$ be the minimum of all the numbers $\dist(S^1, U_i)$ and $\dist(U_i, U_j)$, where $i\not= j$. Set $\varepsilon':= \lambda^{-2}L^{-4}\varepsilon$ and $\rho:=\frac{1}{12L^2}\cdot \min\left\{\varepsilon', \rho_0, \nu\right\}$. 
 Fix $i$ and let $\bar{B}_{i,j} \subset U_i\subset Y_i$, $j=1,\dots, k_i$, be a finite collection of balls obtained from applying Proposition~\ref{prop:covering} to $K_i$, $U_i$, and $Y_i$. Let $A_{i,j}\subset K_i\cap \partial \bar{B}_{i,j}$ be a finite subset which is $2\rho r_{i,j}$-dense in $\partial\bar{B}_{i,j}$, where $r_{i,j}$ denotes the radius of $\bar{B}_{i,j}$. 
 
We now consider the sets $A_{i,j}$ and $\Gamma_{i,j}:= \partial\bar{B}_{i,j}$ as subsets of $(\R^2, |\cdot|)$ and set $A:= S^1\cup\cup_{i,j} A_{i,j}$. Let $\varphi\colon A\to X$ be the restriction of $u$ to $A$ and note that $\varphi$ is $\lambda$-Lipschitz and $\varphi|_{S^1} = c$. Since $A_{i,j}$ is $2L\rho r_{i,j}$-dense in $\Gamma_{i,j}$ and $2L\rho<\nu$ it follows from Lemma~\ref{lem:ext-smooth-Jordan-Lip} and the choice of $\nu$ that the piecewise geodesic extension of $\varphi|_{A_{i,j}}$ to $\Gamma_{i,j}$ is $3\lambda$-Lipschitz. We denote the extended map by $\varphi$ again. From Proposition~\ref{prop:almost-isom} we obtain $$\length(\varphi|_{\Gamma_{i,j}})^2 \leq (1+\varepsilon)^2\cdot \length_{Y_i}(\partial\bar{B}_{i,j})^2 = 4\pi (1+\varepsilon)^2\cdot \hm^2_{Y_i}(\bar{B}_{i,j})$$ and hence 
\begin{equation*}
 \begin{split} 
  \sum_{j=1}^{k_i} \length(\varphi|_{\Gamma_{i,j}})^2 &\leq 4\pi (1+\varepsilon)^2\cdot \hm^2_{Y_i}(U_i) \leq 4\pi (1+\varepsilon)^3\cdot \hm^2_{Y_i}(K_i)\\
  &\leq 4\pi (1+\varepsilon)^{5}\cdot  \Area(u|_{K_i}).
 \end{split}
\end{equation*}
 This shows that $$ \sum_{i=1}^m\sum_{j=1}^{k_i} \length(\varphi|_{\Gamma_{i,j}})^2 \leq 4\pi (1+\varepsilon)^{5}\cdot \left[\Fillarea(c) + \varepsilon\right].$$

Define $\Omega_{i,j}:= B_{i,j}$. Since $\Omega_{i,j}$ is the Jordan domain enclosed by $\Gamma_{i,j}$ we obtain
\begin{equation*}
 \begin{split}
  \left|D \setminus \bigcup_{i,j} \Omega_{i,j}\right| &\leq |E| + |F\setminus \cup_{i=1}^m K_i| + \sum_{i=1}^m \left|K_i\setminus \cup_{j=1}^{k_i} \bar{B}_{i,j}\right|\\
  &< \frac{\varepsilon}{\lambda^2} + \frac{\varepsilon}{\lambda^2}  + 6\varepsilon' L^4\cdot \left|\cup_{i=1}^mK_i\right|\\
  &\leq \frac{\varepsilon}{\lambda^2}\cdot (2 + 6\pi).
 \end{split}
\end{equation*}
 
 The following lemma finishes the proof of Theorem~\ref{thm:discretization}.
 
 \bl
  The map $\varphi$ is $3\lambda$-Lipschitz on $S^1\cup\cup_{i,j}\Gamma_{i,j}$. 
 \el
 
 \begin{proof}
  Let $x\in\Gamma_{i,j}$ and $x'\in \Gamma_{i', j'}$ with $(i,j)\not=(i',j')$. Let $a\in A_{i,j}$ be a nearest point for $x$ and $a'\in A_{i',j'}$ be a nearest point for $x'$. We claim that
 \begin{equation}\label{eq:estimate-dist-Gamma}
  |x-a| + |x'-a'|\leq \frac{1}{3}\cdot |a-a'|.
 \end{equation}
  Indeed, if $i'=i$ then $|x-a| + |x'-a'| \leq 2L\rho(r_{i,j} + r_{i,j'})$ and $$|a-a'|\geq L^{-1}\|a-a'\|_i \geq L^{-1}\varepsilon'\cdot(r_{i,j} + r_{i,j'}),$$ hence \eqref{eq:estimate-dist-Gamma} by the choice of $\rho$. If $i\not=i'$ then $$|x-a| + |x'-a'| \leq 2L\rho(r_{i,j} + r_{i',j'})\leq 4L^2\rho \leq \frac{1}{3}\cdot \rho_0 \leq \frac{1}{3}\cdot |a-a'|.$$ 
 This proves the claim. We note that \eqref{eq:estimate-dist-Gamma} also holds in the case that $x'\in S^1$ and $a'=x'$. Moreover, \eqref{eq:estimate-dist-Gamma} implies that $|a-a'|\leq \frac{3}{2}\cdot|x-x'|$. Since $\varphi$ is $\lambda$-Lipschitz on $A$ and $3\lambda$-Lipschitz on each $\Gamma_{i,j}$ we conclude that 
\begin{equation*}
 \begin{split}
  d(\varphi(x), \varphi(x'))&\leq d(\varphi(x), \varphi(a)) + d(\varphi(a), \varphi(a')) + d(\varphi(a'), \varphi(x'))\\
   &\leq 3\lambda |x-a| + \lambda|a-a'| + 3\lambda|a'-x'|\\
   &\leq 3\lambda|x-x'|.
 \end{split}
\end{equation*}
This completes the proof.
  \end{proof}

\section{Constructing competitors}\label{sec:competitors}

In this section we use Theorem~\ref{thm:discretization} together with Lipschitz extension techniques in order to show that the filling area of a curve in an ultralimit of certain sequences of metric spaces bounds from above the filling areas of approximating curves.

Let $C, L\geq 1$ and $0< r_0\leq\infty$. For each $n\in\N$ let $(X_n, d_n)$ be a proper, geodesic metric space which is $L$-Lipschitz $1$-connected up to some scale and satisfies 
\begin{equation}\label{eq:isop-below-threshhold}
 \delta_{X_n}(r)< C\cdot r^2
\end{equation}
 for all $r\in(0,r_0)$. Suppose furthermore that there exist $\varepsilon_n\in(0,1)$ with $\varepsilon_n\to 0$ as $n\to\infty$ and such that 
\begin{equation}\label{eq:almost-eucl-isop-above-threshhold}
 \delta_{X_n}(r) < \frac{1+\varepsilon_n}{4\pi}\cdot r^2
\end{equation}
 for all $n\in\N$ and all $\varepsilon_n\leq r<r_0$.

Fix a non-principal ultrafilter $\omega$ on $\N$ and a sequence of basepoints $p_n\in X_n$ and denote by $X_\omega$ the ultralimit $(X_n, d_n, p_n)_\omega$. We denote the metric on $X_\omega$ by $d_\omega$. Recall from Section~\ref{subsec:asymp-cones} the definition of bounded sequence of curves and its ultralimit. With the assumptions above we have:

\bt\label{thm:competitors-in-Xn}
  Let $(c_n)$ be a bounded sequence of curves $c_n\colon S^1\to X_n$ with uniformly bounded Lipschitz constants and let $c\colon S^1\to X_\omega$ be given by $c = \lim\nolimits_\omega c_n$. If $\length(c)< r_0$ then for every $\varepsilon>0$ there exists a subset $N\subset\N$ with $\omega(N)= 1$ and such that $$\Fillarealip_{X_n}(c_n) \leq \Fillarea_{X_\omega}(c) + \varepsilon$$ for every $n\in N$.
\et
 
The rest of this section is devoted to the proof of the theorem above. 
By Theorem~\ref{thm:Sobolev-Dehn-stable}, the ultralimit $X_\omega$ satisfies $\delta_{X_\omega}(r)\leq \frac{r^2}{4\pi}$ for all $r\in (0,r_0)$ and hence has property (ET) by Theorem~\ref{thm:eucl-isop-property-ET}. Let $(c_n)$ be a bounded sequence of curves as in the statement of the theorem above and suppose that the ultralimit $c$ satisfies $\length(c)<r_0$. Let $\varepsilon\in(0,1)$. We may assume that $\varepsilon$ is so small that $$(1+\varepsilon)^2\left[\length(c)^2 + \varepsilon\right] <r_0^2$$ and that $\varepsilon<\frac{r_0^2}{16M^2}$, where $M$ is the universal constant appearing after the proof of Lemma~\ref{lem:12-lip} below. 

Let $\lambda\geq 1$ be sufficiently large and, in particular, so large that each $c_n$ is $\lambda$-Lipschitz. By Theorem~\ref{thm:discretization} there exist finitely many pairwise separated smooth convex Jordan curves $\Gamma_1,\dots, \Gamma_k\subset D$ and a $\lambda$-Lipschitz map $\varphi\colon K\to X_\omega$ with $K= S^1\cup \Gamma_1\cup\dots\cup \Gamma_k$ such that $\varphi|_{S^1} = c$ and $$\sum_{i=1}^k \length\left(\varphi|_{\Gamma_i}\right)^2 \leq 4\pi\cdot \Fillarea_{X_\omega}(c) + \varepsilon.$$ Moreover, the Jordan domains $\Omega_i$ enclosed by $\Gamma_i$ satisfy $\left| D\setminus \bigcup_{i=1}^k\Omega_i\right| \leq \frac{\varepsilon}{\lambda^2}$. Set $\Gamma_0:= S^1$ and notice that $$\rho:= \min\left\{\dist(\Gamma_i, \Gamma_j): 0\leq i<j\leq k\right\}>0.$$ Let $\nu>0$ be so small that $10\nu<\rho$ and that Lemma~\ref{lem:ext-smooth-Jordan-Lip} applies with $\nu$ to every $\Gamma_i$. Let $S\subset K$ be a finite set such that $S\cap \Gamma_i$ is $\nu$-dense in $\Gamma_i$ for every $i=0,\dots, k$. For each $n\in\N$ define a map $\varphi_n\colon S\to X_n$ as follows. If $s\in S\cap \Gamma_0$ then set $\varphi_n(s):= c_n(s)$. If $s\in S\setminus \Gamma_0$ then write $\varphi(s)$ as $\varphi(s)= [(x_n)]$ for some bounded sequence of points $x_n\in X_n$ and define $\varphi_n(s):= x_n$ for every $n\in\N$. If $\varphi(s) = \varphi(s')$ for some $s\not=s'$ then we choose the same sequence. 

We will now extend $\varphi_n$ in several steps to a map defined on $\overline{D}$. In each step the map will be called $\varphi_n$. Firstly, extend $\varphi_n$ to $K$ in such a way that $\varphi_n|_{S^1} = c_n$ and such that $\varphi_n|_{\Gamma_i}$ is a piecewise geodesic extension of $\varphi_n|_{S\cap \Gamma_i}$ for $i\geq 1$. Set $$\delta:= \min\{|s-s'|: s,s\in S, s\not= s'\}>0$$ and let $\eta\in(0,1)$ be sufficiently small, to be determined later.
Let $N\subset \N$ be the set of $n\in\N$ such that $C \varepsilon_n^2 \leq \frac{\varepsilon}{k}$ and $\varepsilon_n<\varepsilon$ and $$|d_n(\varphi_n(s),\varphi_n(s')) - d_\omega(\varphi(s), \varphi(s'))| \leq \eta\cdot \delta$$ for all $s,s'\in S$, where $C$ is the constant appearing in \eqref{eq:isop-below-threshhold}. Notice that $\omega(N)=1$. 

\bl\label{lem:12-lip}
 For every $n\in N$ the map $\varphi_n\colon K\to X_n$ is $12\lambda$-Lipschitz.
\el

\begin{proof}
 For distinct points $s,s'\in S$ we have $$d_n(\varphi_n(s), \varphi_n(s')) \leq d_\omega(\varphi(s), \varphi(s')) + \eta\cdot \delta\leq 2\lambda |s-s'|.$$ Thus, by Lemma~\ref{lem:ext-smooth-Jordan-Lip}, the map $\varphi_n|_{\Gamma_i}$ is $6\lambda$-Lipschitz for every $i$. Finally, let $z\in \Gamma_i$ and $z'\in \Gamma_j$ for some $i\not=j$ and let $s\in S\cap \Gamma_i$ and $s'\in S\cap \Gamma_j$ be such that $|z-s|\leq\nu$ and $|z'-s'|\leq\nu$. Since $|s-s'|\leq |z-z'| + 2\nu$ and $|z-z'|\geq \rho>10\nu$ we obtain $$d_n(\varphi_n(z), \varphi_n(z')) \leq 6\lambda\cdot (\nu + |s-s'| +\nu)\leq 6\lambda\cdot(|z-z'| + 4\nu)\leq 12\lambda|z-z'|.$$ This completes the proof.
\end{proof}

Fix $n\in N$. We construct a Lipschitz extension of $\varphi_n$ with suitable area bound as follows. Firstly, by the classical proof of Lipschitz extensions based on Whitney cube decompositions, there exists a countable collection $\mathcal{Q}$ of pairwise almost disjoint closed squares such that $$D\setminus \bigcup_{i=1}^k \overline{\Omega}_i = \bigcup_{Q\in\mathcal{Q}} Q$$ and there exists an $M\lambda$-Lipschitz extension of $\varphi_n$ to the set $K\cup\mathcal{Q}^{(1)}$, see for example the proof of \cite[Theorem 6.4]{LWY16}. Here, $\mathcal{Q}^{(1)}$ denotes the $1$-skeleton of $\mathcal{Q}$ and $M$ is a universal constant. We denote this Lipschitz extension again by $\varphi_n$. We first extend $\varphi_n$ to $\Omega_i$. If $\varphi|_{\Gamma_i}$ is constant then $\varphi_n|_{\Gamma_i}$ is constant by construction and we can extend $\varphi_n$ to a constant map on $\overline{\Omega}_i$. We may thus assume that $\length(\varphi|_{\Gamma_i})>0$ for every $i$. We first note that for every $i\geq 1$ we have $$\length_{X_n}(\varphi_n|_{\Gamma_i}) \leq (1+\varepsilon)\cdot \length(\varphi|_{\Gamma_i}),$$ provided $\eta$ was chosen sufficiently small, depending on $\varepsilon$, $\delta$, the number of points in $S$, and the minimum of the lengths $\length(\varphi|_{\Gamma_j})$. Now fix $i\geq 1$. Since $$\length(\varphi|_{\Gamma_i})^2 \leq 4\pi\cdot \Fillarea_{X_\omega}(c) + \varepsilon \leq \length(c)^2 + \varepsilon$$ we obtain from the above and choice of $\varepsilon$ that $$\length_{X_n}(\varphi_n|_{\Gamma_i})^2 \leq (1+\varepsilon)^2\cdot \length(\varphi|_{\Gamma_i})^2 \leq (1+\varepsilon)^2\cdot[\length(c)^2 + \varepsilon] < r_0^2$$ and so $\length_{X_n}(\varphi_n|_{\Gamma_i})< r_0$. 
 If $\length_{X_n}(\varphi_n|_{\Gamma_i}) \geq \varepsilon_n$ then, by \eqref{eq:almost-eucl-isop-above-threshhold}, there exists a Lipschitz extension of $\varphi_n|_{\Gamma_i}$ to $\overline{\Omega}_i$ with $$\Area_{X_n}(\varphi_n|_{\Omega_i}) \leq \frac{1+\varepsilon}{4\pi}\cdot \length_{X_n}(\varphi_n|_{\Gamma_i})^2 \leq \frac{(1+\varepsilon)^3}{4\pi}\cdot \length(\varphi|_{\Gamma_i})^2.$$ If $\length_{X_n}(\varphi_n|_{\Gamma_i}) < \varepsilon_n$ then, by \eqref{eq:isop-below-threshhold}, there exists a Lipschitz extension of $\varphi_n|_{\Gamma_i}$ to $\overline{\Omega}_i$ with $$\Area_{X_n}(\varphi_n|_{\Omega_i}) \leq C\cdot \length_{X_n}(\varphi_n|_{\Gamma_i})^2 < C\varepsilon_n^2\leq \frac{\varepsilon}{k}.$$ This yields $$\sum_{i=1}^k\Area_{X_n}(\varphi_n|_{\Omega_i})\leq (1+\varepsilon)^3\cdot \Fillarea_{X_\omega}(c) + \varepsilon\cdot \left[1 + (1+\varepsilon)^3\right].$$

Now, we extend $\varphi_n$ to each $Q\in\mathcal{Q}$. Since $X_n$ is $L$-Lipschitz $1$-connected up to some scale it follows that for all $Q$ for which $\diam(Q)$ is sufficiently small (and hence for all but finitely many $Q$), there exists an $M'L\lambda$-Lipschitz extension of $\varphi_n|_{\partial Q}$ to $Q$ and hence $$\Area_{X_n}(\varphi_n|_Q)\leq (M')^2\lambda^2L^2\cdot|Q|,$$ where $M'$ is a universal constant.
For each of the remaining finitely many $Q\in\mathcal{Q}$ we have $|Q|\leq |D\setminus \cup_{i=1}^k \overline{\Omega_i}|\leq \frac{\varepsilon}{\lambda^2}$ and hence $$\length_{X_n}(\varphi_n|_{\partial Q})\leq M\lambda\cdot \length(\partial Q) \leq 4M\sqrt{\varepsilon} <r_0$$ by the choice of $\varepsilon$. The quadratic isoperimetric inequality \eqref{eq:isop-below-threshhold} thus provides a Lipschitz extension of $\varphi_n|_{\partial Q}$ to $Q$ with $$\Area_{X_n}(\varphi_n|_Q)\leq C\cdot \length_{X_n}(\varphi_n|_{\partial Q})^2\leq 16CM^2\lambda^2\cdot |Q|.$$ Consequently, we have $$\sum_{Q\in\mathcal{Q}} \Area_{X_n}(\varphi_n|_Q)\leq M''\lambda^2\cdot \sum_{Q\in\mathcal{Q}} |Q| = M''\lambda^2\cdot \left|D\setminus\bigcup_{i=1}^k\Omega_i\right|\leq M''\varepsilon$$ for a constant $M''$ only depending on $M$, $M'$, $C$, and $L$.

Finally, since the Lipschitz constant of $\varphi_n|_Q$ is uniformly bounded it follows from the classical proof of the Lipschitz extension theorem that the map $\varphi_n$ is Lipschitz continuous on all of $\overline{D}$. Moreover, we conclude from the above that $$\Area_{X_n}(\varphi_n) = \sum_{i=1}^k \Area_{X_n}(\varphi_n|_{\Omega_i}) + \sum_{Q\in\mathcal{Q}} \Area_{X_n}(\varphi_n|_Q) \leq (1+\varepsilon)^3\cdot \Fillarea_{X_\omega}(c) + M''' \varepsilon$$ for a constant $M'''$ only depending on $M$, $M'$, $C$, and $L$. We conclude that $$\Fillarealip_{X_n}(c_n) \leq (1+\varepsilon)^3\cdot \Fillarea_{X_\omega}(c) + M''' \varepsilon$$ for every $n\in N$. Since $\varepsilon>0$ was arbitrary, this completes the proof of Theorem~\ref{thm:competitors-in-Xn}.

\section{Existence of energy and area minimizers in ultralimits}\label{sec:energy-mins}

In this section we solve the classical problem of Plateau in ultralimits of sequences of proper geodesic metric spaces admitting almost Euclidean isoperimetric inequalities. For a fixed proper metric space, the Plateau problem was solved in \cite{LW15-Plateau}. This was extended to a certain class of locally non-compact metric spaces in \cite{GW17-Plateau}. Neither of these results can be applied to the setting we are interested in here.

Given a complete metric space $X$ and a Jordan curve $\Gamma\subset X$ we define $\Lambda(\Gamma, X)$ to be the set of all $v\in W^{1,2}(D, X)$ whose trace has a continuous representative which is a weakly monotone parametrization of $\Gamma$. In other words, $\trace(v)$ has a continuous representative which is the uniform limit of homeomorphisms $c_i\colon S^1\to \Gamma$. The following weak notion of conformality was introduced in \cite{LW15-Plateau}. A map $v\in W^{1,2}(D, X)$ is said to be conformal if for almost every $z\in D$ we have $$\apmd v_z(w) = \apmd v_z(w')$$ for all $w,w'\in S^1$. 

The main result of this section can now be stated as follows.

\bt\label{thm:energy-min-ultralimits}
Let $0<r_0\leq \infty$ and let $(\varepsilon_n)$ be a sequence of positive real numbers tending to $0$. For every $n\in\N$ let $(X_n, d_n)$ be a proper, geodesic metric space satisfying $$\delta_{X_n}(r)\leq \frac{1+\varepsilon_n}{4\pi}\cdot r^2 + \varepsilon_n$$ for all $r\in(0,r_0)$. Let $X_\omega=(X_n, d_n, p_n)_\omega$ be the ultralimit with respect to some non-principal ultrafilter $\omega$ on $\N$ and some sequence of basepoints $p_n\in X_n$. Then for every rectifiable Jordan curve $\Gamma\subset X_\omega$ with $\length(\Gamma)<r_0$ there exists $u\in \Lambda(\Gamma, X_\omega)$ with $$E_+^2(u) = \inf\left\{E_+^2(v): v\in \Lambda(\Gamma, X_\omega)\right\}.$$ Every such $u$ is conformal and also minimizes area among all elements in $\Lambda(\Gamma, X_\omega)$. 
\et

Moreover, $u$ has a unique representative which is locally Lipschitz continuous on $D$ and is continuous on $\overline{D}$. This follows from \cite[Theorem 1.4]{LW15-Plateau} and the fact that $\delta_{X_\omega}(r)\leq \frac{1}{4\pi}\cdot r^2$ for all $r\in (0,r_0)$. 
Notice that unlike in Euclidean space or Riemannian manifolds, in the realm of metric spaces energy minimizers need not be area minimizers, see \cite[Proposition 11.6]{LW15-Plateau}. However, in proper metric spaces satisfying property (ET) energy minimizers are area minimizers as was shown in \cite[Theorem 11.4]{LW15-Plateau}. Of course, $X_\omega$ need not be proper.

The theorem above and the remark after Theorem~\ref{thm:Sobolev-Dehn-stable} yield the following result which, in particular, implies Theorem~\ref{thm:Plateau-asymp-intro}.

\bc\label{cor:energy-min-asymp-cone}
 Let $X$ be a proper, geodesic metric space satisfying  $$\limsup_{r\to\infty} \frac{\delta_X(r)}{r^2} \leq \frac{1}{4\pi},$$  and let
$X_\omega$ be an asymptotic cone of $X$. Then for every rectifiable Jordan curve $\Gamma\subset X_\omega$ there exists $u\in \Lambda(\Gamma, X_\omega)$ with $$E_+^2(u) = \inf\left\{E_+^2(v): v\in \Lambda(\Gamma, X_\omega)\right\}.$$ Every such $u$ is conformal and also minimizes area among all elements in $\Lambda(\Gamma, X_\omega)$. 
\ec

Theorem~\ref{thm:energy-min-ultralimits} can easily be deduced from the following result. Let $X_n$, $X_\omega$, and $\Gamma$ be as in Theorem~\ref{thm:energy-min-ultralimits}.

\bt\label{thm:area-min-ultralimits}
There exists $u\in \Lambda(\Gamma, X_\omega)$ such that $$\Area(u) = \inf\left\{\Area(v): v\in \Lambda(\Gamma, X_\omega)\right\}$$ and such that the image of $u$ is contained in a compact subset of $X_\omega$.
\et

We first provide:

\begin{proof}[Proof of Theorem~\ref{thm:energy-min-ultralimits}]
 Let $u$ be as in Theorem~\ref{thm:area-min-ultralimits} and let $K\subset X_\omega$ be a compact set containing the image of $u$. Then $K$ must contain $\Gamma$. We may thus view $u$ as an element of $\Lambda(\Gamma, K)$. Define a non-empty family of Sobolev maps by $$\Lambda_u:= \left\{v\in \Lambda(\Gamma, K): \Area(v)= \Area(u)\right\}.$$ By the arguments in the proof of \cite[Theorem 7.1]{LW15-Plateau} there exists an element $w\in\Lambda_u$ which minimizes the energy $E_+^2$ among all elements in $\Lambda_u$. By Theorem~\ref{thm:Sobolev-Dehn-stable} we have $\delta_{X_\omega}(r)\leq \frac{r^2}{4\pi}$ for all $r\in (0, r_0)$ and so $X_\omega$ has property (ET) by Theorem~\ref{thm:eucl-isop-property-ET}.  Thus, also $K$ has property (ET). Therefore, \cite[Theorem 11.3]{LW15-Plateau} implies that $w$ is conformal. It follows with \cite[Lemma 7.2]{LW15-Plateau} that $$E_+^2(w) = \Area(w) = \Area(u) \leq \Area(v)\leq E_+^2(v)$$ for every $v\in\Lambda(\Gamma, X_\omega)$. This shows that $w$ is an energy minimizer and an area minimizer in $\Lambda(\Gamma, X_\omega)$.

Finally, let $v$ be any energy minimizer in $\Lambda(\Gamma, X_\omega)$. Then $v$ is conformal by \cite[Theorem 11.3]{LW15-Plateau} and satisfies $$\Area(v) = E_+^2(v) = E_+^2(w) = \Area(w),$$ where $w$ is as above. This shows that $v$ also minimizes area. This completes the proof.
\end{proof}

We turn to the proof of Theorem~\ref{thm:area-min-ultralimits}. Let $r_0$, $\varepsilon_n$, $X_n=(X_n, d_n)$, $p_n$, $\omega$, $X_\omega$, and $\Gamma$ be as in the statement of Theorem~\ref{thm:energy-min-ultralimits}. By Proposition~\ref{prop:good-thickening} we may assume that there exist $C, L\geq 1$ such that $X_n$ is $L$-Lipschitz $1$-connected  up to some scale and satisfies $\delta_{X_n}(r)< Cr^2$ for all $r\in(0,r_0)$ as well as $$\delta_{X_n}(r)< \frac{1+\varepsilon_n}{4\pi}\cdot r^2$$ for all $\varepsilon_n\leq r< r_0$. Here, $(\varepsilon_n)$ is a possibly different sequence but still tends to $0$.

For every $n\in\N$ denote by $\overline{D}_n$ the unit disc $\overline{D}$ equipped with the metric $\varepsilon_n|\cdot|$. Define a metric space $Y_n$ by $Y_n:= X_n\times \overline{D}_n$, where we equip $Y_n$ with the Euclidean product metric, again denoted by $d_n$. Notice that $Y_n$ is proper and geodesic and satisfies $\delta_{Y_n}(r)\leq C'r^2$ for all $r\in(0,r_0)$, where $C'$ only depends on $C$, see \cite[Lemma 3.2]{LW16-harmonic}. View $X_n$ as a subset of $Y_n$ by identifying $X_n$ with $X_n\times\{0\}$. Then the Hausdorff distance between $X_n$ and $Y_n$ tends to zero as $n\to\infty$. Hence, $X_\omega$ is isometric to the ultralimit $(Y_n, d_n, \bar{p}_n)_\omega$, where the basepoints $\bar{p}_n\in Y_n$ are defined by $\bar{p}_n:=(p_n, 0)$.

Let $c\colon S^1\to \Gamma$ be a constant speed parametrization of the rectifiable Jordan curve $\Gamma\subset X_\omega$. Recall that $\length(\Gamma)<r_0$ by assumption. By \cite[Corollary 2.6]{LWY16} there exists a bounded sequence $(c_n)$ of curves $c_n\colon S^1\to X_n$ with uniformly bounded Lipschitz constants such that $c = \lim\nolimits_\omega c_n$ and $$\length_{X_n}(c_n) + 2\pi \varepsilon_n <r_0$$ for all sufficiently large $n$. For all such $n$ define an injective curve in $Y_n$ by $$\bar{c}_n(z):= (c_n(z), z)$$ for every $z\in S^1$ and notice that $(\bar{c}_n)$ is a bounded sequence with uniformly bounded Lipschitz constants. Moreover, $c = \lim\nolimits_\omega \bar{c}_n$ and $\length_{Y_n}(\bar{c}_n) \leq \length_{X_n}(c_n) + 2\pi \varepsilon_n < r_0$. 

Denote by $\Gamma_n$ the Jordan curve in $Y_n$ given by $\Gamma_n:= \bar{c}_n(S^1)$. By \cite{LW15-Plateau} there exists $u_n\in \Lambda(\Gamma_n, Y_n)$ which is continuous on $\overline{D}$, minimizes area in $\Lambda(\Gamma_n, Y_n)$ and minimizes energy among all area minimizers. In particular, it follows from \cite[Theorem 6.2]{LW15-Plateau} that $u_n$ satisfies $$E_+^2(u_n) \leq 2\cdot \Area(u_n) \leq 2C'\cdot \length_{Y_n}(\bar{c}_n)^2,$$ which is uniformly bounded. Fix distinct points $q_1, q_2, q_3\in S^1$. After possibly composing with a conformal diffeomorphism of $D$ we may assume that $u_n$ satisfies the $3$-point condition $u_n(q_i) = \bar{c}_n(q_i)$ for $i=1,2,3$. By the proof of \cite[Proposition 8.7]{LW15-Plateau} there exists for every $s\in (0,1)$ some $L_s>0$ such that $u_n$ is $(L_s, \alpha)$-H\"older continuous on $\overline{B}(0, s)$ for every $n\in\N$, where $\alpha = \frac{1}{8\pi C'}$.

\bl\label{lem:unif-cpt-An}
 The sequence of metric spaces $(A_n, d_n)$, where $A_n = u_n(\overline{D})$, is uniformly compact in the sense of Gromov.
\el

\begin{proof}
 By \cite[Section 1]{LW-intrinsic}, the set $A_n$ is the image under a $1$-Lipschitz map of a geodesic metric space $Z_n$ which is homeomorphic to $\overline{D}$ and satisfies the following properties. Firstly, the Hausdorff $2$-measure of $Z_n$ equals $\Area(u_n)$ and the length of the boundary circle $\partial Z_n$ equals $\length_{Y_n}(\bar{c}_n)$. Notice that both these quantities are bounded from above by some number $M$ which does not depend on $n$. Secondly, for all $z\in Z_n$ and $0\leq r\leq \dist(z,\partial Z_n)$ we have $$\hm_{Z_n}^2(B(z,r))\geq C''r^2$$ for a constant $C''$ only depending on $C'$. In particular, the diameter of $Z_n$ is bounded from above by $M + 2\sqrt{M/C''}$. Moreover, for every $k\in \N$, there exists some $\frac{M}{k}$-dense subset of $Z_n$ which has at most $C'''k^2$ elements, where $C'''$ only depends on $C'$, see \cite[Corollary 8.10]{LW-intrinsic}. From this the statement of the lemma follows.
\end{proof}

By Gromov's compactness theorem for metric spaces \cite{Gro81-poly} there exists a compact metric space $(Z, d_Z)$ and isometric embeddings $\varphi_n\colon A_n\hookrightarrow Z$ for all $n\in\N$. Define continuous maps $v_n\colon \overline{D}\to Z$ by $v_n:= \varphi_n\circ u_n$. Let $v$ be the ultralimit of the sequence $(v_n)$, thus $$v(z):= \lim\nolimits_\omega v_n(z)$$ for every $z\in \overline{D}$. Notice that $v$ is $(L_s, \alpha)$-H\"older continuous on $\overline{B}(0,s)$ for every $s\in(0,1)$. Define injective Lipschitz curves by $\gamma_n:= \varphi_n\circ \bar{c}_n$ and let $\gamma$ be the ultralimit of $(\gamma_n)$. It follows that $\gamma$ is a Lipschitz curve in $Z$.

\bl
Define a subset $A\subset Z$ by $A:= \left\{ \lim\nolimits_\omega \varphi_n(a_n): a_n\in A_n\right\}$. Then the map $\psi\colon A\to X_\omega$ given by $$\psi(\lim\nolimits_\omega \varphi_n(a_n)):= [(a_n)]$$ is well-defined and an isometric embedding.
\el

\begin{proof}
 Let $(a_n)$ be a sequence with $a_n\in A_n$ for all $n$. Then $$\sup_{n\in\N} d_n(a_n, \bar{p}_n)<\infty$$ because the diameter of $A_n$ is uniformly bounded by Lemma~\ref{lem:unif-cpt-An} and $(\bar{c}_n)$ is a bounded sequence. Now, if $a = \lim\nolimits_\omega\varphi_n(a_n)$ and $a'= \lim\nolimits_\omega\varphi_n(a'_n)$ are two points in $A$ then 
 \begin{equation*}
   d_\omega([(a_n)], [(a'_n)]) = \lim\nolimits_\omega d_n(a_n, a'_n) = \lim\nolimits_\omega d_Z(\varphi_n(a_n), \varphi_n(a'_n)) = d_Z(a, a'),
 \end{equation*}
which shows that $\psi$ is well-defined and an isometric embedding.
\end{proof}

Since $X_\omega$ is a complete metric space, $\psi$ extends to an isometric embedding from the closure $\overline{A}$ of $A$ to $X_\omega$. We denote this map by $\psi$ again. We notice that $\psi\circ\gamma = c$ and hence $\gamma$ is injective. We denote by $\Gamma'$ the image of $\gamma$, which is thus a rectifiable Jordan curve in $Z$.

\bl\label{lem:properties-v-in-Z}
 The map $v$ belongs to $\Lambda(\Gamma', Z)$ and satisfies $\Area(v)\leq \Fillarea_{X_\omega}(c)$.
\el

\begin{proof}
We first show that $v\in W^{1,2}(D, Z)$ with $\Area(v)\leq \Fillarea_{X_\omega}(c)$. For this, let $\varepsilon>0$ and let $N\subset \N$ be a subset with $\omega(N)=1$ as in Theorem~\ref{thm:competitors-in-Xn}, when applied to $X_n$ and $c_n$. Since the curves $c_n$ and $\bar{c}_n$ can be connected by a Lipschitz annulus in $Y_n$ of area at most proportional to $\varepsilon_n$ and $\varepsilon_n\to 0$ it follows that for all but finitely many $n\in N$ we have $$\Area_{Y_n}(u_n) \leq \Fillarea_{X_n}(c_n) + \varepsilon \leq \Fillarea_{X_\omega}(c) + 2\varepsilon.$$ Let $S\subset \overline{D}$ a countable dense set such that $S\cap S^1$ is dense in $S^1$. By the definition of ultralimit there exists a strictly increasing sequence of numbers $n_k\in N$ such that $v_{n_k}(s) \to v(s)$ for all $s\in S$ and $\gamma_{n_k}(s) \to \gamma(s)$ for all $s\in S\cap S^1$. It follows that $\gamma_{n_k}$ converges to $\gamma$ uniformly on $S^1$ and that $v_{n_k}$ converges to $v$ locally uniformly on $D$ and, in particular, the convergence is in $L^2(D, Z)$. Since $v_{n_k}\in W^{1,2}(D, Z)$ and $E_+^2(v_{n_k}) = E_+^2(u_{n_k})$ is uniformly bounded it follows from \cite[Theorem 1.13]{KS93} that $v\in W^{1,2}(D, Z)$ and from \cite[Corollary 5.8]{LW15-Plateau} that $$\Area_Z(v)\leq \liminf_{k\to\infty} \Area_Z(v_{n_k}) \leq \Fillarea_{X_\omega}(c) + 2\varepsilon.$$ Since $\varepsilon>0$ was arbitrary we see that $\Area_Z(v)\leq \Fillarea_{X_\omega}(c)$.

It remains to show that $\trace(v)$ is a weakly monotone parametrization of $\Gamma'$. Firstly, it follows from Lemma~\ref{lem:gen-equi-cpt} below that the family $\{v_{n_k}|_{S^1}: k\in\N\}$ is equi-continuous. Thus, after possibly passing to a further subsequence, we may assume that $v_{n_k}|_{S^1}$ converges uniformly to a weakly monotone parametrization $\gamma'$ of $\Gamma'$. Since $v_{n_k}|_{S^1}$ converges in $L^2(S^1, Z)$ to $\trace(v)$ by \cite[Theorem 1.12.2]{KS93} it follows that $\trace(v) = \gamma'$. This proves that $v\in \Lambda(\Gamma', Z)$ and completes the proof.
\end{proof}

Notice that the image of $v$ lies in the compact set $\overline{A}$. Hence, the map $u:= \psi\circ v$ belongs to $\Lambda(\Gamma, X_\omega)$, has image in the compact set $K:= \psi(\overline{A})$ containing $\Gamma$, and satisfies $$\Area_{X_\omega}(u)\leq \Fillarea_{X_\omega}(c).$$ By \cite[Lemma 4.8]{LW-intrinsic}, we have $\Fillarea_{X_\omega}(c)\leq \Area_{X_\omega}(w')$ for all $w'\in\Lambda(\Gamma, X_\omega)$, which shows that $u$ minimizes area among all elements in $\Lambda(\Gamma, X_\omega)$. This completes the proof of Theorem~\ref{thm:area-min-ultralimits}.

The following slight generalization of \cite[Proposition 7.4]{LW15-Plateau} was used in the proof of Lemma~\ref{lem:properties-v-in-Z}.

\bl\label{lem:gen-equi-cpt}
Let $Z$ be a complete metric space. Let $\gamma_k\colon S^1\to Z$ be continuous, injective curves converging uniformly to an injective curve $\gamma\colon S^1\to Z$. Set $\Gamma_k:= \gamma_k(S^1)$ and let $M> 0$. Let $q_1,q_2,q_3\in S^1$ be distinct points and suppose $v_k\in \Lambda(\Gamma_k, Z)$ satisfies the $3$-point condition $\trace(v_k)(q_i) = \gamma_k(q_i)$ for $i=1,2,3$ and $E_+^2(v_k)\leq M$ for all $k$. Then the family $\{\trace(v_k): k\in\N\}$ is equi-continuous.
\el

 In particular, a subsequence of $(\trace(v_k))$ converges uniformly to a weakly monotone parametrization of $\Gamma = \gamma(S^1)$.
 
\begin{proof}
 We first notice that for every $\varepsilon>0$ there exists $\delta>0$ such that if $k\in\N$ and $x,y\in \Gamma_k$ satisfy $d(x,y)<\delta$ then one of the two segments of $\Gamma_k$ between $x$ and $y$ lies in the ball $B(x,\varepsilon)$. This together with the Courant-Lebesgue lemma and the $3$-point condition now implies that the family $\{\trace(v_k): k\in\N\}$ is equi-continuous, exactly as in the case of a single Jordan curve.
\end{proof}

\section{The main result and its consequences}\label{sec:asymp-CAT}

The following may be considered the main result of this paper. It generalizes Theorem~\ref{thm:asymp-CAT-intro} stated in the introduction and will also be used to prove Theorem~\ref{thm:intro-bigger-const-CAT} and has other consequences.

\bt\label{thm:main-CAT}
 Let $0<r_0\leq \infty$ and let $(\varepsilon_n)$ be a sequence of positive real numbers tending to $0$. For every $n\in\N$ let $(X_n, d_n)$ be a proper, geodesic metric space satisfying $$\delta_{X_n}(r)\leq \frac{1+\varepsilon_n}{4\pi}\cdot r^2 + \varepsilon_n$$ for all $r\in(0,r_0)$. Let $X_\omega=(X_n, d_n, p_n)_\omega$ be the ultralimit with respect to some non-principal ultrafilter $\omega$ on $\N$ and some sequence of basepoints $p_n\in X_n$. Then every geodesic triangle in $X_\omega$ of perimeter strictly smaller than $r_0$ is ${\rm CAT}(0)$.
\et

Theorem~\ref{thm:main-CAT} together with the remark after Theorem~\ref{thm:Sobolev-Dehn-stable} implies Theorem~\ref{thm:asymp-CAT-intro}. The proof of Theorem~\ref{thm:main-CAT} is done by combining Theorem~\ref{thm:energy-min-ultralimits} with the arguments in the proof of the main result in \cite{LW-isoperimetric}. We will actually use a strengthening of one of the main theorems in \cite{LW-isoperimetric} established in \cite{LW-param}.

\begin{proof}
 By Theorem~\ref{thm:Sobolev-Dehn-stable} we have that $$\delta_{X_\omega}(r)\leq \frac{1}{4\pi}\cdot r^2$$ for every $r\in(0,r_0)$. Let $\Gamma\subset X_\omega$ be a geodesic triangle of perimeter strictly smaller than $r_0$. We want to show that $\Gamma$ is ${\rm CAT}(0)$. We may assume that $\Gamma$ defines a Jordan curve in $X_\omega$, see the proof of \cite[Lemma 3.1]{LW-isoperimetric}. 
 
 By Theorem~\ref{thm:energy-min-ultralimits}, there exists $u\in \Lambda(\Gamma, X_\omega)$ which minimizes the Reshetnyak energy $E_+^2$ and the area among all elements of $\Lambda(\Gamma, X_\omega)$ and which is conformal. Moreover, $u$ has a representative which is continuous on $\overline{D}$ by \cite[Theorem 1.4]{LW15-Plateau}.  Thus, by \cite[Section 1]{LW-intrinsic}, there exists a geodesic metric space $Z$, called the intrinsic minimal disc associated with $u$, and a $1$-Lipschitz map $\bar{u}\colon Z\to X_\omega$ with the following properties. Firstly, the space $Z$ is homeomorphic to $\overline{D}$ and the restriction of $\bar{u}$ to the boundary circle $\partial Z$ is an arc-length preserving homeomorphism from $\partial Z$ onto $\Gamma$. In particular, $\length_Z(\partial Z) = \length_{X_\omega}(\Gamma)<r_0$. Secondly, $\hm^2_Z(Z)= \Area_{X_\omega}(u)$ and every Jordan domain $\Omega\subset Z$ satisfies $$\hm^2_Z(\Omega) \leq \frac{1}{4\pi}\cdot \length(\partial \Omega)^2.$$ Notice that \cite[Theorem 1.2]{LW-intrinsic} only asserts this inequality for Jordan domains $\Omega\subset Z$ with $\length_Z(\partial \Omega)<r_0$. However, in the above this also holds when $\length_Z(\partial \Omega)\geq r_0$ because in this case $$\hm_Z^2(\Omega)\leq \hm^2_Z(Z) = \Area_{X_\omega}(u)\leq \frac{1}{4\pi}\cdot \length_{X_\omega}(\Gamma)^2< \frac{1}{4\pi}\cdot r_0^2 \leq \frac{1}{4\pi}\cdot \length_Z(\partial \Omega)^2.$$ Now, it follows from \cite[Corollary 1.5]{LW-param} that $Z$ is a ${\rm CAT}(0)$-space. The proof of \cite[Lemma 3.3]{LW-isoperimetric} shows that $\Gamma$ is ${\rm CAT}(0)$. This completes the proof. 
\end{proof}

The proof of Theorem~\ref{thm:intro-bigger-const-CAT} follows from Theorem~\ref{thm:main-CAT} together with Proposition~\ref{prop:ultralimit-CAT-consequence-sequence}:

\begin{proof}[Proof of Theorem~\ref{thm:intro-bigger-const-CAT}]
 Suppose by contradiction that the statement is wrong. Then there exist $\nu\in(0,1)$, a sequence $(r_n)$ of positive real numbers, and a sequence of proper, geodesic metric spaces $(X_n, d_n)$ with the following property. For each $n\in\N$ the space $X_n$ satisfies $$\deltalip_{X_n}(r)\leq \frac{1+\frac{1}{n}}{4\pi}\cdot r^2$$ for all $r\in(0,r_n)$ but $X_n$ contains a geodesic triangle $\Delta_n$ of perimeter $s_n< (1-\nu)r_n$ for which the ${\rm CAT}(0,\nu\cdot s_n)$-condition fails. For each $n\in\N$ define a new metric by $\bar{d}_n:= s_n^{-1}d_n$ and define the rescaled metric space $Y_n:= (X_n, \bar{d}_n)$. Now, view $\Delta_n$ as a triangle in $Y_n$. Its perimeter in $Y_n$ is $1$ and it fails the ${\rm CAT}(0, \nu)$-condition in $Y_n$. Hence, by Proposition~\ref{prop:ultralimit-CAT-consequence-sequence}, there is a non-principal ultrafilter $\omega$ on $\N$ and a sequence of basepoints $p_n\in Y_n$ such that the ultralimit $Y_\omega=(Y_n, \bar{d}_n, p_n)_\omega$ contains a geodesic triangle of perimeter at most $1$ which fails to be ${\rm CAT}(0)$. However, this contradicts Theorem~\ref{thm:main-CAT} since each $Y_n$ satisfies
$$\deltalip_{Y_n}(r)\leq \frac{1+\frac{1}{n}}{4\pi}\cdot r^2$$ for all $r\in(0,\frac{r_n}{s_n})$ and $\frac{r_n}{s_n} > \frac{1}{1-\nu}>1$. This concludes the proof.
\end{proof}

Notice that the Lipschitz Dehn function in Theorem~\ref{thm:intro-bigger-const-CAT} and its proof can be replaced by the Sobolev Dehn function. Theorem~\ref{thm:main-CAT} also implies the following result which can be regarded as a coarse analog of Theorem~\ref{thm:intro-bigger-const-CAT}. The proof is very similar to the one above.

\bt\label{thm:intro-coarse-Dehn-CAT}
 For all $r_0>0$ and $\nu\in(0,1)$ there exists $\varepsilon>0$ with the following property. If $X$ is a proper, geodesic metric space satisfying $$\delta_X(r)\leq \frac{1}{4\pi}\cdot r^2 + \varepsilon$$ for all $r\in(0,r_0)$ then every geodesic triangle in $X$ of perimeter at most $(1-\nu)r_0$ is ${\rm CAT}(0, \nu)$. 
\et

\begin{proof}
 Suppose by contradiction that the statement is wrong. Then there exist $r_0>0$, $\nu\in(0,1)$, and a sequence of proper, geodesic metric spaces $(X_n, d_n)$ with the following property. Each $X_n$ satisfies $$\delta_{X_n}(r)\leq \frac{1}{4\pi}\cdot r^2 + \frac{1}{n}$$ for all $r\in (0,r_0)$ but the ${\rm CAT}(0,\nu)$-condition fails for some geodesic triangle in $X_n$ of perimeter at most $(1-\nu)r_0$. Proposition~\ref{prop:ultralimit-CAT-consequence-sequence} thus implies that for some non-principal ultrafilter $\omega$ on $\N$ and some sequence of basepoints $p_n\in X_n$ the ultralimit $X_\omega=(X_n, d_n, p_n)_\omega$ must contain a geodesic triangle of perimeter at most $(1-\nu)r_0$ which fails to be ${\rm CAT}(0)$. However, this contradicts Theorem~\ref{thm:main-CAT} and finishes the proof.
\end{proof}

The following proposition shows that the constant $\frac{1}{4\pi}$ in Theorem~\ref{thm:asymp-CAT-intro} is optimal.

\bp
 For every $\varepsilon>0$ there exist some $2$-dimensional non-Euclidean normed space $X$ and $\frac{1}{4\pi}<C<\frac{1}{4\pi}+\varepsilon$ such that $\deltalip_X(r) = Cr^2$ for all $r\geq 0$.
\ep

Notice that as a non-Euclidean normed space, $X$ is not ${\rm CAT}(0)$, see \cite[Proposition II.1.14]{BrH99}.

\begin{proof}
Every $2$-dimensional normed space $X$ satisfies 
\begin{equation}\label{eq:Dehn-normed-plane}
 \deltalip_X(r) = Cr^2
\end{equation}
 for some constant $C\geq \frac{1}{4\pi}$ and for all $r\geq 0$, with $C=\frac{1}{4\pi}$ if and only if $X$ is Euclidean. This follows from inequality \eqref{eq:isop-non-Euclidean} and the area formula. Thus, choosing a non-Euclidean norm $\|\cdot\|$ on $\R^2$ which is sufficiently close to the standard Euclidean one we obtain that $X= (\R^2, \|\cdot\|)$ satisfies \eqref{eq:Dehn-normed-plane} with a constant $C$ which is arbitrarily close to and strictly bigger than $\frac{1}{4\pi}$. 
\end{proof}

We end this paper with:

\begin{proof}[Proof of Theorem~\ref{thm:no-converse}]
Consider the pinwheel tilling of the Euclidean plane $\R^2$ by isometric triangles of side lengths  $1$, $2$, $\sqrt{5}$ constructed in \cite{Rad94-pin}. Notice that each triangle has area equal to $1$. Let $G\subset\R^2$ be the graph consisting of the edges of the triangles in the pinwheel tiling and equip $G$ with the length metric which we denote by $d_G$. Let $X$ be the geodesic metric space obtained by gluing spherical caps onto (the boundaries of) the triangles in $G$. Then $X$ is biLipschitz homeomorphic to $\R^2$ and contains $G=(G, d_G)$ isometrically. Since $G$ is at finite Hausdorff distance from $X$ it follows that the asymptotic cones of $X$ and $G$ are isometric. Moreover, \cite[Theorem 2]{RS96} shows that for every $\varepsilon>0$ there exists $R>0$ such that $$|x-y|\leq d_G(x,y) \leq (1+\varepsilon)\cdot |x-y|$$ whenever $x,y\in G$ satisfy $|x-y|\geq R$. From this it follows that the Euclidean plane $\R^2$ is the unique asymptotic cone of $G$ and thus also of $X$.

It remains to show that $X$ satisfies \eqref{eq:no-converse-ineq}. For this, let $\varepsilon>0$ be suitably small, to be determined below. It follows from \cite[Theorem 1]{RS96} that for every sufficiently large $r>0$ there exists a Jordan curve $\Gamma\subset G$ whose length satisfies $$(1-\varepsilon)\cdot r \leq \length(\Gamma) \leq r$$ and such that $\Gamma$ encloses at least $\left(\frac{1}{4\pi} - \varepsilon\right) \cdot \length(\Gamma)^2$ triangles of the pinwheel tiling. Since each spherical cap in $X$ has Hausdorff $2$-measure bigger than $3$ it follows that, as a subset of $X$, the curve $\Gamma$ encloses a Jordan domain $\Omega\subset X$ of Hausdorff measure at least $$\hm^2_X(\Omega) \geq 3\cdot \left(\frac{1}{4\pi} - \varepsilon\right) \cdot \length(\Gamma)^2 \geq 3\cdot \left(\frac{1}{4\pi} - \varepsilon\right) \cdot (1-\varepsilon)^2 \cdot r^2.$$ Thus, if $\varepsilon>0$ was chosen sufficiently small then $$\delta_X(r) = \deltalip_X(r) \geq \frac{1}{2\pi} \cdot r^2$$ for all $r>0$ large enough. This proves \eqref{eq:no-converse-ineq} and completes the proof.
\end{proof}

\def\cprime{$'$} \def\cprime{$'$} \def\cprime{$'$}

\end{document}